\documentclass[10pt]{article} 

\usepackage{amsmath}
\usepackage{amssymb} 
\usepackage{latexsym} 
\usepackage{epic,eepic,epsfig,psfrag} 
\usepackage{amscd}  
\usepackage{index}  

\newcommand{\QED}{\hspace*{\fill}$\Box$\medskip} 
\def\one{\hbox{1\hskip-2.7pt l}}

\def\rd{{\rm d}}
\def\rT{{\rm T}}

\def\rG{{\rm G}}
\def\p{\phi}

\def\a{\alpha}
\def\b{\beta} 
\def\d{\delta} 

\def\ep{\varepsilon} 
\def\e{\eta} 
\def\g{\gamma}
\def\th{\theta} 
\def\r{\rho}
 
\def\t{\tau} 
 
\def\x{\xi} 
 
\def\n{\nu} 
\def\m{\mu} 
 
\def\o{\omega} 
\def\l{\lambda} 
\def\D{\Delta} 
 
\def\G{\Gamma} 

\def\S{\Sigma}
\def\Si{\Sigma} 
\def\Th{\Theta} 
\def\X{\Xi} 
\def\Om{\Omega} 
\def\P{\Phi} 
\def\cg{\mathfrak{g}}
\def\cA{{\mathcal A}}
 
\def\cC{{\mathcal C}}

\def\cG{{\mathcal G}}

\def\cL{{\mathcal L}} 
 
\def\cm{{\mathfrak m}}

\def\cR{{\mathcal R}}

\def\CS{{\mathcal C}{\mathcal S}}
\def\E{{\mathcal E}}
\def\R{{\mathbb R}} 
  
\def\N{{\mathbb N}} 
 
\def\H{{\mathbb H}}

\def\half{{\textstyle{\frac 12}}} 

\def\laplace{\Delta} 
\def\Pr{{\bf Proof:}\;} 
\def\st{\: \big| \:}

\def\dt{{\rm d}t}
\def\ds{{\rm d}s}

\def\dy{{\rm d}y}
\def\dph{{\rm d}\p}
\def\dth{{\rm d}\th}
\def\dr{{\rm d}r}
\def\pd{\partial}
\def\comp{\circ}

\def\ta{{\tilde\alpha}}
\def\tg{{\tilde g}}
\def\tu{{\tilde u}}
\def\tw{{\tilde w}}

\def\tx{{\tilde\x}}
\def\tA{{\tilde A}}

\def\tR{{\tilde R}}

\def\tX{{\tilde\X}}

\def\br{{\bar r}}
\def\bph{{\bar\p}}
\def\bp{{\bar p}}
\def\bq{{\bar q}}

\def\bm{{\bar\m}}
\def\bep{{\bar\ep}}
\def\bt{{\bar\t}}
\def\la{\langle\,}
\def\ra{\,\rangle}
\def\dvol{\,\rd{\rm vol}}

\def\HF{{\rm HF}}

\def\SU{{\rm SU}}

\def\tint{\textstyle\int}
\def\hol{{\rm hol}}

\newtheorem{dfn}{Definition}[section] 
\newtheorem{lem}[dfn]{Lemma} 
\newtheorem{prp}[dfn]{Proposition} 
\newtheorem{thm}[dfn]{Theorem} 
\newtheorem{rmk}[dfn]{Remark}

\begin{document}

\bibliographystyle{plain}

\author{Katrin Wehrheim \\
\vspace{-1.5mm}
{\small wehrheim@princeton.edu; (609)2584234;}\\
\vspace{-2mm}
{\small Princeton University, Fine Hall, Princeton NJ 08544-1000}
}

\title{Anti-self-dual instantons with Lagrangian boundary conditions II: Bubbling}

\maketitle

\begin{abstract}
We study bubbling phenomena of anti-self-dual instantons on $\H^2\times\S$,
where $\S$ is a closed Riemann surface. 
The restriction of the instanton to each boundary slice $\{z\}\times\S$, $z\in\pd\H^2$ 
is required to lie in a Lagrangian submanifold of the moduli space of flat connections over $\S$
that arises from the restrictions to the boundary of flat connections on a handle body.

We establish an energy quantization result for sequences of 
instantons with bounded energy near $\{0\}\times\S$: Either their curvature is in fact uniformly bounded in 
a neighbourhood of that slice (leading to a compactness result) 
or there is a concentration of some minimum quantum of energy.
We moreover obtain a removable singularity result for instantons with finite energy in a punctured 
neighbourhood of $\{0\}\times\S$. 
This completes the analytic foundations for the construction of an instanton Floer homology for
$3$-manifolds with boundary. 
This Floer homology is an intermediate object in the program proposed by Salamon for the proof of the Atiyah-Floer 
conjecture for homology-$3$-spheres.

In the interior case, for anti-self-instantons on $\R^2\times\S$, our methods provide a new approach to the 
removable singularity theorem by Sibner-Sibner for codimension $2$ singularities with a holonomy condition.
\end{abstract}

\section{Introduction}

The aim of this paper is to complete the analytic foundations for the definition
of instanton Floer homology groups $\HF^{\rm inst}_*(M,L_Y)$ begun in \cite{W elliptic}.
Here $M$ is a compact, oriented $3$-manifold with boundary $\pd M=\Si$ and $L_Y\subset M_\Si$ 
is a (singular) Lagrangian submanifold of the moduli space $M_\Si$ of flat connections 
on the trivial $\SU(2)$-bundle over $\Si$.
This Lagrangian is $L_Y:=\cL_Y/\cG^{1,p}(\Si)$, where \hbox{$\cL_Y\subset\cA^{0,p}(\Si)$} is a 
Lagrangian submanifold of the space of $L^p$-connections given by the $L^p$-closure of the flat 
connections on a handle body $Y$ restricted to $\pd Y=\Si$.
Here $\cG^{1,p}(\Si)$ is the group of $W^{1,p}$-gauge transformations, and $L_Y$ is actually 
independent of the choice of $p>2$. 

This Floer homology serves as intermediate object in the program proposed by Salamon
\cite{Sa1} for the proof of the Atiyah-Floer conjecture for homology-$3$-spheres.

Fukaya \cite{Fu} was the first to suggest the use of Lagrangian boundary conditions in order 
to define a Floer homology for $3$-manifolds with boundary. His setup uses nontrivial
bundles (where the moduli spaces of flat connections are smooth manifolds) and thus cannot 
immediately be used in the context of the Atiyah-Floer conjecture, where the bundles are
necessarily trivial and thus the moduli spaces of flat connections are singular.
Our approach is to define $\HF^*_{\rm inst}(M,L_Y)$ from the moduli spaces of anti-self-dual instantons
on $\R\times M$ with Lagrangian boundary condition in $L_Y$, i.e.\ from the gauge equivalence classes
of connections $\X\in\cA(\R\times M)$ satisfying the boundary value problem
\begin{equation}\label{intro bvp}
\left\{\begin{array}{l}
F_\X + *F_\X = 0,\\
\X|_{\{s\}\times\pd M} \in\cL_Y \quad\forall s\in\R .
\end{array}\right.
\end{equation}
Note that the boundary condition is nonlocal: It firstly asserts the local condition that the connection
is flat on each boundary slice; but secondly its holonomy has to vanish around those loops in 
$\Si$ that are contractible in $Y$, which is a global condition.

In \cite{W elliptic} we describe this approach in full detail and we establish the elliptic 
theory for this boundary value problem (allowing for a larger class of Lagrangian boundary 
conditions). Fix $p>2$, then every $W^{1,p}_{\rm loc}$-solution is gauge equivalent to a smooth 
solution and the following analogue of Uhlenbeck compactness is true: 
Every sequence of solutions with locally $L^p$-bounded curvature is gauge equivalent to a 
sequence that contains a $\cC^\infty$-convergent subsequence.

In this paper we address the question of bubbling: What happens if a sequence of solutions
has bounded energy $\int_{\R\times M} |F_\X|^2 <\infty$ but its curvature $F_\X$ is not locally
$L^p$-bounded for any $p>2$? 

In the case of a $4$-manifold without boundary this question is answered by the compactification
of the moduli space of anti-self-dual instantons leading to the Donaldson invariants of 
smooth $4$-manifolds~\cite{D2} and to the instanton Floer homology groups of closed
$3$-manifolds~\cite{F1}. This compactification is described in terms of trees of anti-self-dual
instantons on $S^4$ that 'bubble off' at isolated points on the original $4$-manifold.
In the case of the present boundary value problem, we do not attempt this compactification but
only establish the relevant facts for the definition of the Floer homology groups.
So the answer comes in two parts (that will be stated more precisely in theorems~\ref{thm A} and 
\ref{thm B}):\\

\noindent
{\bf Energy quantization:} If the curvature is not uniformly bounded near an interior point 
$x\in\R\times{\rm int}\,M$ or near a boundary slice $\{s\}\times\Si\subset\R\times\pd M$, then
there is a minimum energy $\ep_0>0$ that concentrates at this point or slice.\\

\noindent
{\bf Removal of singularities:} Every smooth finite energy solution on the complement of an 
interior point or a boundary slice can be put into a gauge in which it extends to a solution
over the full manifold.\\

In the case of interior points, these are the two wellknown analytic ingredients for the 
compactification of the moduli space (see e.g.\ \cite{U1} for Uhlenbeck's removable singularity
theorem). The anti-self-dual instantons on $S^4$ are obtained by rescaling the connections near 
the bubbling point $x$. The limit object then is an instanton on $\R^4$ whose singularity at 
infinity can be removed resulting in an instanton on a nontrivial bundle over $S^4$.

In the case of bubbling at the boundary, one might also find instantons on $S^4$ bubbling off
at boundary points. These would arise from sequences of solutions $\X^\n$ and interior points $x^\n$ 
with distance $t^\n\to 0$ to the boundary $\R\times\pd M$, where the curvature 
$|F_{\X^\n}(x^\n)|=(R^\n)^2$ blows up at a rate such that $R^\n t^\n\to\infty$.
If $R^\n t^\n$ stays bounded, then the standard rescaling construction will lead to 
anti-self-dual instantons on increasingly large domains of the half space. 
In \cite{Sa1} it was conjectured that there is an energy quantization for the limit objects -- 
anti-self-dual instantons on the half space.
However, the local rescaling construction looses the global part of the boundary condition.
With only the slicewise flatness as boundary condition, one cannot expect to obtain better 
convergence than weak $W^{1,p}$-convergence (for any $p<\infty$) up to the boundary.
In the interior, one of course has smooth convergence, and thus might find a nontrivial limit
object. However, in case $R^\n t^\n\to 0$, even the limit object might be trivial if the 
blowup is in the curvature part for which one does not have $\cC^0$-convergence up to the 
boundary.
\footnote{
Writing $\X=\P\ds+\Psi\dt+A$ near the boundary $\{t=0\}$ and assuming $p>4$, one obtains
$W^{2,p}$-bounds for $\X$ except for the second $\pd_s,\pd_t$-derivatives of the connections 
$A(s,t)$ on the $\Si$-slices. 
These bounds suffice to obtain $\cC^0$-convergence for the curvature component $F_A$, but not 
for $B_s=\pd_s A + \rd_A\P$.
The latter requires full $W^{2,p}$-bounds, which would only result from a Lagrangian boundary 
condition coupled with the Cauchy-Riemann equation for $A$ as a function with values in 
$\cA^p(\S)$, c.f.\ \cite{W elliptic}.
}

This discussion suggests a more global analysis of the bubbling phenomenon taking into account
the full $\Si$-slices and localizing only in the two other variables. An adapted rescaling
construction seems to lead to holomorphic discs in the space of connections over $\Si$ (with
the Hodge operator as complex structure) with Lagrangian boundary conditions. We do not have
a precise convergence statement. However, we were able to prove the corresponding energy quantization 
result by purely analytic means -- after all using partial convergence results for the naive 
local rescaling construction described above.\\

Before giving the precise statements of our main results we introduce the setup and some
basic notation. (For more details on gauge theory and the notation used here see \cite{W Cauchy} or \cite{W}.)
Throughout this paper, we are working in a small neighbourhood of a boundary slice of a
Riemannian $4$-manifold with a boundary space-time splitting in the sense of 
\cite[Def~1.2]{W elliptic}. So we are considering the following local model.

We denote by $B_r(x_0)\subset\R^n$ the closed ball of radius $r>0$ centered at 
$x_0\in\R^n$.
The intersection of a ball with the half space 
$$
\H^n:=\{(s_1,\ldots,s_{n-1},t)\in\R^n \;|\; t\geq 0 \} 
$$
is denoted by
$$ 
D_r(x_0):= B_r(x_0)\cap\H^n .
$$
Moreover, we write $D:=D_{r_0}(0)\subset\H^2$ for the $2$-dimensional 
half ball centered at $0$ of some fixed radius $r_0$.
Next, let $\Si$ be a closed Riemann surface. Now the local model is the trivial 
$\SU(2)$-bundle over the Riemannian $4$-manifold
$$
(\, D\times\Si \,,\, \ds^2 + \dt^2 + g_{s,t} \,).
$$
Here $g_{s,t}$ is a family of metrics on $\Si$ that varies smoothly with
$(s,t)\in D$.
We will call any metric of this type a {\bf metric of normal type}.

For all purposes in this paper, we can replace $\SU(2)$ by a general compact, connected, 
and simply connected Lie group $\rG$.
Now a $\rG$-connection on $D\times\S$ is a $1$-form \hbox{$\X\in\Om^1(D\times\Si,\cg)$}
with values in the Lie algebra $\cg$. We will write $\cA(X)$ for the space of smooth 
connections over a manifold $X$, then $\cA_{\rm flat}(X)$ denotes the space of smooth 
flat connections, 
and $\cG(X)$ is the space of smooth gauge transformations on $X$ 
(i.e.\ maps to $\rG$).
The Sobolev spaces of connections and gauge transformations are denoted by
\begin{align*}
\cA^{k,p}(X) &= W^{k,p}(X,\rT^*X\otimes\cg), \\
\cG^{k,p}(X) &= W^{k,p}(X,\rG).
\end{align*}
We will be dealing with anti-self-dual instantons on $D\times\S$ that satisfy a 
Lagrangian boundary condition as follows.
Let $p>2$ and fix a handle body $Y$ with boundary $\pd Y=\Si$, then the following
Lagrangian submanifold is introduced in \cite[Lemma~4.6]{W Cauchy},
$$
\cL_Y \,:=\; {\rm cl}_{L^p}\, \bigl\{ A\in\cA_{\rm flat}(\S) \st \exists 
                              \tA\in\cA_{\rm flat}(Y) : \tA|_\S=A \bigr\} 
\;\subset\;\cA^{0,p}(\Si).
$$
We consider the following boundary value problem for connections
$\X\in\cA(D\times\Si)$
\begin{equation}\label{bvp}
\left\{\begin{array}{l}
F_\X + *F_\X = 0,\\
\X|_{(s,0)\times\Si} \in\cL_Y \quad\forall s\in [-r_0,r_0] .
\end{array}\right.
\end{equation}
The compactness result \cite[Thm~B]{W elliptic} for this boundary value problem
can be phrased as follows for the local model. 
Here ${\rm int}(D)={\rm int}(B_{r_0}(0))\cap\H^2$ 
denotes the interior in the topology of $\H^2$.

\begin{thm} {\bf (Compactness) \cite{W elliptic}}  \label{thm cp} \\
Let $p>2$ and let $g^\n$ be a $\cC^\infty$-convergent sequence of metrics of normal type on $D\times\S$.
Suppose that $\X^\n\in \cA(D\times\S)$ is a sequence of solutions of (\ref{bvp}) with respect to the
metrics $g^\n$ such that $\|F_{\X^\n}\|_{L^p(D\times\S)}$ is uniformly bounded.

Then there exists a subsequence (again denoted by $\X^\n$) and a sequence of gauge transformations 
$u^\n\in\cG(D\times\S)$ such that $u^{\n\;*}\X^\n$ converges uniformly with all derivatives on every
compact subset of ${\rm int}(D)\times\S$.
\end{thm}

Next, we state the energy quantization result that will be proven in section~\ref{quantization}.

\begin{thm} {\bf (Energy quantization)}  \label{thm A} \\
Let $r_0>0$ and let $\cm$ be a $\cC^\infty$-compact set of metrics of normal type on 
$D\times\Si$. Then there exists a constant $\ep_0>0$ such that the following holds.

Let $\X^\n\in\cA(D\times\Si)$ be a sequence of solutions of (\ref{bvp}) with 
respect to metrics $g^\n\in\cm$. 
Suppose that for all $\d>0$
$$
\sup_\n \; \bigl\| F_{\X^\n} \bigr\|_{L^\infty(D_\d(0)\times\Si)} \;=\; \infty . 
$$
Then after taking a subsequence there exist $(s^\n,t^\n)\to 0$ and $\ep^\n\to 0$ such 
that
$$
\int_{D_{\ep^\n}(s^\n,t^\n)\times\Si} \bigl| F_{\X^\n} \bigr|^2 \;>\; \ep_0 .
$$
\end{thm}

\begin{rmk} \hspace{1mm} \label{rmk thm A} \\
\vspace{-5mm} 
\begin{enumerate}
\item 
By theorem~\ref{thm cp} the assumptions in theorem~\ref{thm A} imply that for a subsequence and with any $p>2$
one has for all $\d>0$
$$
\sup_\n \; \bigl\| F_{\X^\n} \bigr\|_{L^p(D_\d(0)\times\Si)} \;=\; \infty  .
$$
\item
With the stronger assumption in (i) it suffices to consider a $\cC^3$-compact set of metrics in the
theorem, as will be seen in the proof.
By following through the proof of theorem~\ref{thm cp},
in particular \cite[Thm~2.6]{W elliptic}, one can moreover
check that the set of metrics in theorem~\ref{thm A} only needs to be $\cC^5$-compact.
\end{enumerate}
\end{rmk}

To see (i) note that otherwise one would find a sequence $\X^\n$ of solutions with respect to a 
$\cC^\infty$-convergent sequence of metrics $g^\n$ and constants $C$, \hbox{$\d>0$} such that 
\hbox{$\|F_{\X^\n}\|_{L^p(D_{2\d}\times\S)}\leq C$} but
\hbox{$\|F_{\X^\n}\|_{L^\infty(D_\d\times\S)}\to\infty$}.
Due to the $L^p$-bounded curvature one would then find a subsequence and 
gauges in which the connections converge uniformly on 
$D_\d\times\S$. Since the norm of the curvature is gauge invariant, this contradicts the above divergence.
In fact, we will need to make the stronger assumption in (i) for some $2<p<3$ in order to deduce the energy
quantization directly. (This is why we had to establish theorem~\ref{thm cp} in \cite{W elliptic} in the technically more 
difficult case $2<p\leq 4$.)

With this stronger assumption the structure of the proof of theorem~\ref{thm A} will be similar to an argument in the interior 
case, where it is possible to obtain the energy quantization result independently of the removal of singularities 
and of any geometric knowledge about energies of instantons on $S^4$. 
This argument just uses a wellknown mean value inequality for the Laplace operator and will also be explained in 
section \ref{quantization}.
In our case we will need a mean value inequality up to the boundary at which we cannot simply reflect the function.
Instead, we will use a mean value inequality for functions with a control on the Laplacian and on the
normal derivative at the boundary, which we introduce in \cite{W mean}. 
The following result from \cite{W mean} should give an idea of this type of a priori estimate -- 
in the actual proof, we will need a slightly different, more special version.

\begin{lem} \label{lem a}
For every $n\geq 2$ there exists a constant $C$ such that for all \hbox{$A,B\geq 0$} there exists 
$\m(A,B)>0$ with the following significance.

Let $D_r(y)\subset\H^n$ be the Euclidean $n$-ball in the half space of radius $r>0$ and center 
$y\in\H^n$. 
Suppose that $e\in\cC^2(D_r(y),[0,\infty))$ satisfies
\begin{equation*}
\left\{\begin{array}{ll}
\laplace e &\leq  B e^{\frac {n+2}n} , \\
\frac\pd{\pd\n}\bigr|_{\pd\H^n} e \hspace{-2mm} &\leq A e^{\frac {n+1}n} ,
\end{array}\right.
\qquad\text{and}\qquad
\int_{D_r(y)} e < \m(A,B) .
\end{equation*}
Then
$$
e(y) \leq C r^{-n} \int_{D_r(y)} e .
$$
\end{lem}

With the energy quantization established, every sequence of solutions of (\ref{intro bvp}) with bounded energy
converges smoothly on the complement of finitely many interior points and boundary slices
(modulo gauge and taking a subsequence). Now the remaining key analytic point for the 
definition of the Floer homology groups is to show that the limit object -- after gauge --
gives rise to a new solution, that will have less energy.
At the interior points, this is Uhlenbeck's removable singularity theorem \cite[Thm~4.1]{U1}.
For the boundary slices, this requires the following removal of codimension-$2$-singularities 
that will be proven in section~\ref{singularity}.
Here again $D\subset\H^2$ denotes the standard closed half ball with center $0$ and some 
fixed radius $r_0>0$, and we introduce the punctured half balls
$$
D^*_r := D_r(0) \setminus\{0\}, \qquad\qquad  D^*:= D^*_{r_0} = D\setminus \{0\}.
$$

\begin{thm} {\bf (Removal of singularities for boundary slices)} \label{thm B} \\
Let $\X\in\cA(D^*\times\S)$ be a smooth connection with finite energy 
$\int_{D^*\times\S}|F_\X|^2<\infty$ and suppose that it satisfies
\[
\left\{\begin{array}{l}
*F_\X + F_\X = 0,\\
\X|_{(s,0)\times\Si} \in\cL_Y \quad\forall s\in[-r_0,0)\cup (0,r_0] .
\end{array}\right.
\]
Then there exists a gauge transformation $u\in\cG(D^*\times\S)$ such that $u^*\X$ extends to a smooth
connection and solution of (\ref{bvp}) on $D\times\S$.
\end{thm}

Both the energy quantization and the removal of singularities rely on the specific form 
of the Lagrangian boundary condition: Connections in $\cL_Y\subset\cA^{0,p}(\Si)$ are 
extended from $\pd Y=\Si$ to flat connections on $Y$ with the $L^2$-norm on $\Si$ controlling 
the $L^3$-norm on $Y$.
The corresponding linear and nonlinear extension results are given in the following lemma 
and are proven in section~\ref{extension}.

\begin{lem} \hspace{1mm}  \label{lem b}
There exists a constant $C_Y$ such that the following holds.
\begin{enumerate}
\item 
For every smooth path $A : (-\ep,\ep)\to \cL_Y\cap\cA(\S)$ there exists another path
$\tA : (-\ep,\ep)\to \cA_{\rm flat}(Y)$ with $\pd_s\tA(0)|_{\pd Y}=\pd_s A(0)$ 
such that
$$
\|\pd_s\tA(0)\|_{L^3(Y)} \leq C_Y \|\pd_s A(0)\|_{L^2(\S)} .
$$
\item
For all $A_0, A_1\in\cL_Y\cap\cA(\S)$ there exist
$\tA_0,\tA_1\in\cA_{\rm flat}(Y)$ with $A_i=\tA_i|_{\pd Y}$ 
such that
\begin{equation} \label{b est}
\|\tA_0-\tA_1\|_{L^3(Y)} \leq C_Y \|A_0-A_1\|_{L^2(\S)} .
\end{equation}
\end{enumerate}
\end{lem}

\pagebreak

\begin{rmk} \label{rmk b}
The constant $C_Y$ in lemma~\ref{lem b} can be chosen uniform
for a $\cC^0$-neighbourhood of metrics on $Y$ and the induced metrics on $\S=\pd Y$.
\end{rmk}

This can be seen by using a fixed metric for the construction of the extensions.
The $L^2(\S)$- and $L^3(Y)$-norms for different metrics are then equivalent with a 
small factor for $\cC^0$-close metrics.

The nonlinear extension in (ii) allows to define a local Chern-Simons
functional for short arcs from $\cL_Y$ to $\cL_Y$: 
We consider smooth paths \hbox{$A:[0,\pi]\to\cA(\S)$} with endpoints $A(0),A(\pi)\in\cL_Y$ .
For such paths lemma~\ref{lem b}~(ii) provides extensions $\tA(0),\tA(\pi)\in\cA_{\rm flat}(Y)$ of 
$A(0),A(\pi)$ that satisfy (\ref{b est}). We pick any such extensions to define
\begin{align}
\CS(A) &:= - \half\int_0^\pi \int_\S \la A \wedge \pd_\p A \ra \;\dph  \label{CS int}\\
&\quad  + \tfrac 1{12} \int_Y \la \tA(0) \wedge [\tA(0)\wedge\tA(0)] \ra 
                         - \la \tA(\pi) \wedge [\tA(\pi)\wedge\tA(\pi)]  \ra  .  \nonumber
\end{align}
Here the notations $[ \cdot \cdot ]$ and $\la \cdot \cdot \ra$ indicate that the values of the 
differential forms are paired via the Lie bracket and an equivariant inner product on $\cg$ 
respectively. 
This is the actual Chern-Simons functional on 
$\bar Y \cup_{\{0\}\times\S} [0,\pi]\times\S \cup_{\{\pi\}\times\S} Y$ 
of the connection given by $\tA(0)$, $A$, and $\tA(\pi)$ on the 
different parts. (Here $\bar Y$ denotes $Y$ with the reversed orientation.)
The extensions $\tA(0)$ and $\tA(\pi)$ could both vary by
gauge transformations that are trivial on $\pd Y=\S$.
So the connection on the above closed manifold might also vary by a gauge 
transformation (that is trivial on the middle part).
The Chern-Simons functional however does not vary under gauge transformations
that are homotopic to $\one$, and it only changes by multiples of $4\pi^2$ for
others.\footnote
{
This constant is correct for $\rG=\SU(2)$ with $\la\x,\e\ra={\rm tr}(\x^*\e)$.
For a general Lie group we can achieve the same constant by scaling the inner product 
appropriately.
}
In fact, if we restrict to short paths, then we will see in section~\ref{isoperimetric} 
that our local Chern-Simons functional is welldefined and satisfies an isoperimetric
inequality.

\begin{lem} {\bf (Isoperimetric inequality)} \label{lem c} \\
There exists $\ep>0$ such that for all smooth paths $A:[0,\pi]\to\cA(\S)$ with 
$A(0),A(\pi)\in\cL_Y$ and $\int_0^\pi \|\pd_\p A\|_{L^2(\S)}\leq\ep$ 
the local Chern-Simons functional (\ref{CS int}) is welldefined and
satisfies
$$
| \CS(A) |
\;\leq\; \left( \int_0^\pi \bigl\| \pd_\p A \bigr\|_{L^2(\S)} \,\dph \right)^2 .
$$
\end{lem}

The significance of the local Chern-Simons functional for theorem~\ref{thm B} is in
the fact that the energy of the connection can be expressed by this functional.
The isoperimetric inequality will then provide a control on the rate of 
decay of the energy on small neighbourhoods of the singularity.
%
%
%
%
%
%
%
%
%
This can be combined with mean value inequalities as in lemma~\ref{lem a} 
to obtain estimates on the connection (in a specific gauge) near the singularity. 
Finally, we will be able to remove the singularity using a cutoff construction
and the compactness result, theorem~\ref{thm cp}.

Note that in our approach all bubbling at the boundary is treated globally, even if it could be
described as an instanton on $S^4$ bubbling off at the boundary.
In fact, the energy quantization result also holds for interior slices (i.e.\  
$\{s\}\times\{t\}\times\S \subset \R\times{\rm int}\,M$ in a tubular neighbourhood 
$\R\times[0,\ep)\times\S$ of $\R\times\pd M$).
This description of the bubbling phenomena would then require a removable singularity result for 
anti-self-dual instantons with a singularity of codimension $2$. An obviously necessary condition for
this result is that the limit holonomy around the singularity vanishes almost everywhere.
It was shown by Sibner-Sibner \cite[Thm~5.2]{Si} and Rade \cite[Thm~2.1]{Ra} 
that this condition is in fact sufficient.
Moreover, the fact that interior bubbling only occurs at isolated points shows that the holonomy 
condition is satisfied at interior slices. This is of little use in our context, so we stick to a
pointwise description of interior bubbling.

However, our techniques for the removal of slice singularities at the boundary also give rise to an
alternative approach to the Sibner-Sibner result for interior slices.
In fact, this approach might lead to a general normal form in terms of the limit holonomy
for finite energy anti-self-dual instantons with a singularity of codimension 2.
(This question was raised by Kronheimer and Mrowka in \cite{KM}.)
However, in this paper, we only consider a special case in which we obtain a largely simplified
proof of the removal of singularities. 
This proof is given in section~\ref{singularity}. In order to state the result we
denote by $B$ the standard closed ball with center $0$ and some fixed radius $r_0>0$, 
and we introduce the punctured ball $B^*$,
$$
B:=B_{r_0}(0)\subset\R^2 , \qquad\qquad
B^*:= B\setminus \{0\} .
$$
Introducing polar coordinates $(r,\p)$ on $B^*$ one can write
any connection on $B^*\times\S$ in the form $\X=R\dr + \P \rd\p + A$, 
where $A$ is a family of $1$-forms on $\S$.
The holonomy condition in $\cite{Si}$ is equivalent to the existence of a gauge in which
\begin{equation*} \label{eq:holonomy condition}
\int_0^{2\pi} \bigl\| \P(r,\p) \bigr\|_{L^2(\S)}^2 \rd\p  \;\;\underset{r\to 0}{\longrightarrow}\; 0 .
\end{equation*}
We will make the stronger assumption that in fact there is a gauge in a neighbourhood of the
singular slice in which $\P\equiv 0$.


\begin{rmk} {\bf (Removal of singularities for interior slices) \cite{Si,Ra}} \label{rmk C} \\
Let $\X\in\cA(B^*\times\S)$ be a smooth anti-self-dual connection with 
finite energy $\int_{B^*\times\S}|F_\X|^2<\infty$
and suppose that $\X$ is gauge equivalent to a connection on $B^*\times\S$ with 
$\P\equiv 0$.
Then there exists a gauge transformation $u\in\cG(B^*\times\S)$ such that $u^*\X$ 
extends to a smooth anti-self-dual connection on $B\times\S$.
\end{rmk}

\subsubsection*{Acknowledgements}

Dietmar Salamon has contributed more than he cares to claim to these results
-- a lot of expertise and both encouragement and criticism.
Fengbo Hang filled a gap by explaining theorem~\ref{thm p ext} to me.
This research was partially supported by the Swiss National Science Foundation.

\pagebreak

\section{Energy quantization}
\label{quantization}

The energy quantization result for anti-self-dual instantons 
at interior points could be phrased as follows (in the special case of a Euclidean
metric).

\begin{thm} \label{thm quant0}
There exists a constant $\ep_0>0$ such that the following holds.

Let $B:=B_{r_0}(0)\subset\R^4$ be the Euclidean $4$-ball of radius $r_0>0$ and let 
$\X^\n\in\cA(B)$ be a sequence of anti-self-dual connections. 
Suppose that 
$$
\sup_\n \;\bigl\| F_{\X^\n} \bigr\|_{L^\infty(B_\d(0))} \;=\; \infty 
\qquad\quad \forall \d>0 .
$$
Then after taking a subsequence there exist $B\ni x^\n \to 0$ and $\ep^\n\to 0$ such 
that for all $\n\in\N$
$$
\int_{B_{\ep^\n}(x^\n)} \bigl| F_{\X^\n} \bigr|^2 \;>\; \ep_0 .
$$
\end{thm}

This is of course a wellknown result in gauge theory. Here we give a purely analytic 
proof that does not use the removable singularity result. This exhibits a general method
for establishing energy quantization whenever one has a (nonlinear) bound on the Laplacian
of the energy density, and this implies a mean value inequality on balls of small energy.
In our case, this mean value inequality will be provided by the following wellknown result
(see e.g.\ \cite{W mean}).

\begin{prp}
\label{prp mean int}
For every $n\in\N$ there exist constants $C$, $\m>0$, and $\d>0$ such that the 
following holds.

Let $\R^n$ be equipped with a metric $g$ such that $\|g-{\emph \one}\|_{W^{1,\infty}}\leq\d$.
Let $B_r(0)\subset\R^n$ be the geodesic ball of radius $0<r\leq 1$. 
Suppose that $e\in\cC^2(B_r(0),[0,\infty))$ satisfies for some $A,B\geq 0$
$$
\laplace e \leq A e + B e^{\frac {n+2}n}
\qquad\text{and}\qquad
\int_{B_r(0)} e < \m B^{-\frac n2} .
$$
Then
$$
e(0) \leq C \bigl( A^{\frac n2} + r^{-n} \bigr) \int_{B_r(0)} e .
$$
\end{prp}

\noindent
{\bf Proof of theorem \ref{thm quant0}:}\;
By assumption one can find a subsequence and points $B\ni x^\n \to 0$ such that 
$R^\n:=| F_{\X^\n}(x^\n)|^{\frac 12} \to \infty$.
We pick a sequence $\ep^\n\to 0$ such that still $\ep^\n R^\n \to \infty$.
Now consider the energy density functions \hbox{$e^\n= | F_{\X^\n} |^2 : B \to [0,\infty)$}.
One can check (see (\ref{BoWei}) below) that $\laplace e^\n \leq 8 (e^\n)^{\frac 32}$.
Let $\m>0$ be the constant from the mean value inequality proposition~\ref{prp mean int}, 
then the theorem holds with $\ep_0=\frac\m{64}$. 
Indeed, for all sufficiently large $\n\in\N$
(such that $B_{\ep^\n}(x^\n)\subset B$) we either have
$\int_{B_{\ep^\n}(x^\n)} e^\n > \ep_0$, or by means of proposition~\ref{prp mean int}
$$
(R^\n)^4 \;=\;e^\n(x^\n) \;\leq\; C (\ep^\n)^{-4}\int_{B_{\ep^\n}(x^\n)} e^\n
$$
and thus 
$(\ep^\n R^\n)^4 \leq C \ep_0 .$
Since $\ep^\n R^\n \to \infty$ the latter can only be true for finitely many $\n\in\N$.
\QED

The proof of theorem~\ref{thm A} will run along similar lines. Here the mean value inequality
(with a boundary condition) will be applied to the functions 
$\|F_{\X^\n}\|_{L^2(\Si)}^2$ that are defined on $D=D_{r_0}(0)\subset\H^2$.
So firstly, we need to show that the assumption in theorem~\ref{thm A}, i.e.\ no local 
uniform bound for the curvature near the slice $\{0\}\times\Si$, 
actually implies a blowup of the above function (the slicewise $L^2$-norm of the curvature) 
at $0\in\H^2$. 
Here remark \ref{rmk thm A}~(i) is crucial: It asserts that in fact there is no local 
$L^p$-bound for the curvature near $\{0\}\times\Si$ for any $p>2$.
From this stronger assumption (we need $p<3$), lemma~\ref{lem L2 bound} below will then 
imply the blowup of $\|F_{\X^\n}\|_{L^2(\Si)}^2$.

The underlying analytic facts of this lemma and the whole proof of theorem~\ref{thm A} will be mean value 
inequalities for both $\|F_{\X^\n}\|_{L^2(\Si)}^2$ (on a 2-dimensional domain with
boundary) and $|F_{\X^\n}|^2$ (on a 4-dimensional domain).
So we shall first calculate the Laplacians and normal derivatives of these functions.
For that purpose we write the connection in the splitting
$$
\X = A + \P\ds + \Psi\dt ,
$$
where $A:D\to\Om^1(\Si,\cg)$ and $\P,\Psi:D\to\Om^0(\Si,\cg)$.\footnote
{
Note that this notation differs from \cite{W elliptic}, where we wrote
$A=B+\P\ds+\Psi\dt$.
}
By $\rd_A$ and $\rd_A^*$ we then denote the families (parametrized by $(s,t)\in D$) 
of operators on $\S$ corresponding to $A(s,t)$.
Moreover, we introduce the covariant derivatives 
$$ 
\nabla_s := \pd_s + [\P,\cdot], \qquad
\nabla_t := \pd_t + [\Psi,\cdot].
$$ 
Now the components of the curvature are $F_A$ and
\begin{align*}
B_s :=\;\pd_s A -\rd_A\P
   &\;=\; [\nabla_s,\rd_A],\\
B_t :=\;\pd_t A -\rd_A\Psi
   &\;=\; [\nabla_t,\rd_A],\\
\pd_t\P -\pd_s\Psi + [\Psi,\P] &\;=\; [\nabla_t,\nabla_s] .
\end{align*}
The Bianchi identity $\rd_\X F_\X = 0$ becomes in this splitting 
\begin{align*}
\nabla_s F_A = \rd_A B_s, \qquad
\nabla_t F_A = \rd_A B_t, \qquad
\nabla_s B_t - \nabla_t B_s = \rd_A [\nabla_t,\nabla_s], 
\end{align*}
and the anti-self-duality equation is
\begin{align*}
*B_s = B_t , \qquad
*F_A = [\nabla_t,\nabla_s] .
\end{align*}

\begin{lem} \label{lem Laplacians}
There is a constant $C$ (varying continuously with the metric of normal type in
the $\cC^2$-topology) such that for all solutions $\X\in\cA(D\times\Si)$ of (\ref{bvp})
\begin{align*}
\laplace \bigl| F_\X \bigr|^2
&\;\leq\; C \bigl| F_\X \bigr|^2 + 8 \bigl| F_\X \bigr|^3 ,   \phantom{\int_\S}\\
\laplace \bigl\| F_\X \bigr\|_{L^2(\Si)}^2 
&\;\leq\; C \bigl\| F_\X \bigr\|_{L^2(\Si)}^2 
   - 20 \la F_A \,,\, [B_s\wedge B_s] \ra_{L^2(\Si)}  \\
&\;\leq\; C \bigl( 1 + \bigl\| F_A \bigr\|_{L^\infty(\S)} \bigr) \bigl\| F_\X \bigr\|_{L^2(\Si)}^2  , \\
-\tfrac\pd{\pd t}\bigr|_{t=0} \bigl\| F_\X \bigr\|_{L^2(\Si)}^2 
&\;\leq\; C \bigl\| B_s \bigr\|_{L^2(\Si)}^2 
        - 4\int_\Si \la \nabla_s B_s \wedge B_s \ra  \\
&\;\leq\; C \bigl( \bigl\| B_s \bigr\|_{L^2(\Si)}^2 + \bigl\| B_s \bigr\|_{L^2(\Si)}^3 \bigr) .
\end{align*}
\end{lem}
\Pr
The anti-self-duality equation together with the Bianchi identity
gives
\begin{align*}
\nabla_s B_s + \nabla_t B_t 
&= * \bigl( - \nabla_s B_t + \nabla_t B_s \bigr) - (\pd_s*) B_t + (\pd_t*) B_s \\ 
&= - * \rd_A * F_A  - (\pd_s*) B_t + (\pd_t*) B_s .
\end{align*}
Using this identity we obtain
\begin{align*}
&\bigl( \nabla_s^2 + \nabla_t^2 \bigr) B_s \\
&= \nabla_s \bigl( -\nabla_t B_t  - * \rd_A * F_A  - (\pd_s*) B_t + (\pd_t*) B_s \bigr)
 + \nabla_t \bigl( \nabla_s B_t  - \rd_A * F_A  \bigr) \\ 
&= [*F_A , B_t] 
 - *\rd_A * \nabla_s F_A - *[B_s,*F_A] 
 -  \rd_A * \nabla_t F_A -  [B_t,*F_A] \\
&\quad  
 -(\pd_s*) \bigl(\rd_A *F_A + \nabla_s B_t \bigr)
 +(\pd_t*) \nabla_s B_s
 - *\rd_A (\pd_s*)F_A - \rd_A (\pd_t*) F_A \\
&\quad 
 - (\pd_s^2 *) B_t + (\pd_s\pd_t*) B_s \\
&= \rd_A^* \rd_A B_s + \rd_A \rd_A^* B_s - 3 *[B_s,*F_A] 
 - (\pd_s^2 *) B_t + (\pd_s\pd_t*) B_s \\
&\quad 
 -(\pd_s*) \nabla_t B_s +(\pd_t*) \nabla_s B_s
 - *\rd_A (\pd_s*)F_A - \rd_A (\pd_t*) F_A ,
\end{align*}
\begin{align*}
\bigl( \nabla_s^2 + \nabla_t^2 \bigr) F_A
&= \nabla_s \rd_A B_s + \nabla_t \rd_A B_t \\
&= \rd_A \bigl( \nabla_s B_s + \nabla_t B_t \bigr) + [B_s\wedge B_s] + [B_t\wedge B_t] \qquad\qquad\qquad\\
&= \rd_A \rd_A^* F_A + 2 [B_s\wedge B_s].
\end{align*}
Continuing these calculations leads to the Bochner-Weitzenb\"ock formula
(c.f.\ \cite[Thm~3.10]{BL}) for anti-self-dual connections
$$
0 \;=\; \bigl(\rd_\X\rd_\X^* + \rd_\X^*\rd_\X\bigr) F_\X 
\;=\; \nabla_\X^*\nabla_\X F_\X + F_\X \comp ({\rm Ric}\wedge g + 2R) + \cR^\X(F_\X).
$$
The quadratic term $\cR^\X(F_\X)\in\Om^2(D\times\Si,\cg)$ can be expressed 
with the help of a local orthonormal frame $(e_1,\ldots,e_4)$ of $\rT(D\times\Si)$ as
$$
\cR^\X(F_\X) (X,Y) = 
\sum_{j=1}^4 \bigl\{ [F_\X(e_j,X) , F_\X(e_j,Y)] - [F_\X(e_j,Y) , F_\X(e_j,X)] \bigr\}.
$$
This gives the first estimate
\begin{align}
\laplace \bigl| F_\X \bigr|^2
&= -2 \bigl| \nabla_\X F_\X \bigr|^2 + 2 \la F_\X \,,\, \nabla_\X^*\nabla_\X F_\X \ra 
\nonumber \\
&\leq -2 \la F_\X \,,\, F_\X \comp ({\rm Ric}\wedge g + 2R) \ra 
      -2 \la F_\X \,,\, \cR^\X(F_\X) \ra   \label{BoWei} \\
&\leq C \bigl| F_\X \bigr|^2 + 8 \bigl| F_\X \bigr|^3 .  \nonumber
\end{align}
Here the constant $C$ depends on the Ricci transform ${\rm Ric}$ and the
scalar curvature $R$ of the metric $g$. 
It can thus be chosen uniform for a $\cC^2$-neighbourhood of the fixed metric.

The purpose of the calculations in the beginning is the following identity:
\begin{align*}
- \tfrac 14 \laplace \bigl\| F_\X \bigr\|_{L^2(\Si)}^2 
&= \bigl\| \nabla_s F_A \bigr\|_{L^2}^2 
 + \bigl\| \nabla_t F_A \bigr\|_{L^2}^2 
 + \bigl\| \nabla_s B_s \bigr\|_{L^2}^2 
 + \bigl\| \nabla_t B_s \bigr\|_{L^2}^2 \\ 
&\quad
 + \la F_A \,,\, \bigr( \nabla_s^2 + \nabla_t^2 \bigl) F_A \ra_{L^2(\Si)} 
 + \la B_s \,,\, \bigr( \nabla_s^2 + \nabla_t^2 \bigl) B_s \ra_{L^2(\Si)} \\
&\quad
 + \la *F_A \,,\, (\pd_s^2*) F_A \ra_{L^2(\Si)} 
 + \la *B_s \,,\, (\pd_s^2*) B_s \ra_{L^2(\Si)} \\
&\quad
 + \la (\pd_s*)F_A \,,\, * \nabla_s F_A \ra_{L^2(\Si)} 
 + \la (\pd_s*)B_s \,,\, * \nabla_s B_s \ra_{L^2(\Si)} 
\end{align*}
\begin{align*}
&= \bigl\| \nabla_s B_s \bigr\|_{L^2(\Si)}^2 
 + \bigl\| \nabla_t B_s \bigr\|_{L^2(\Si)}^2 
 + \bigl\| \rd_A B_s \bigr\|_{L^2(\Si)}^2 
 + \bigl\| \rd_A^* B_s \bigr\|_{L^2(\Si)}^2 \\ 
&\quad
 + \bigl\| \nabla_s F_A \bigr\|_{L^2(\Si)}^2 
 + \bigl\| \nabla_t F_A \bigr\|_{L^2(\Si)}^2 
 + \bigl\| \rd_A^* F_A \bigr\|_{L^2(\Si)}^2 
 + 5 \la F_A \,,\, [B_s\wedge B_s] \ra_{L^2(\Si)} \\
&\quad
 - \la B_s \,,\, *(\pd_s^2 *) B_s + (\pd_s^2 *)*B_s - (\pd_s\pd_t*) B_s \ra_{L^2(\Si)}
 + \la *F_A \,,\, (\pd_s^2*) F_A \ra_{L^2(\Si)} \\
&\quad
 + \la (\pd_s*) B_s \,,\, \nabla_t B_s + * \nabla_s B_s \ra_{L^2(\Si)}
 - \la (\pd_t*) B_s \,,\, \nabla_s B_s \ra_{L^2(\Si)} \\
&\quad
 + 2 \la \rd_A B_s \,,\,  * (\pd_s*)F_A  \ra_{L^2(\Si)}
 - \la \rd_A^* B_s \,,\, (\pd_t*) F_A \ra_{L^2(\Si)} .
\end{align*}
This yields the second inequality
$$
\laplace \bigl\| F_\X \bigr\|_{L^2(\Si)}^2 \\
\leq C \bigl(\bigl\| B_s \bigr\|_{L^2(\Si)}^2 + \bigl\| F_A \bigr\|_{L^2(\Si)}^2 \bigr)
  - 20 \la F_A \,,\, [B_s\wedge B_s] \ra_{L^2(\Si)} .
$$
Here the constant $C$ depends on the second derivatives of $g_{s,t}$ and its inverse.
Using the Bianchi identity, the anti-self-duality equation, and in addition
the boundary condition $F_A\bigr|_{t=0}=0$ we obtain for the normal derivative as 
claimed
\begin{align*}
- \tfrac 14 \tfrac\pd{\pd t}\bigr|_{t=0} \bigl\| F_\X \bigr\|_{L^2(\Si)}^2 
&= - \bigl( \la F_A \,,\, \nabla_t F_A \ra_{L^2(\Si)}
 + \la B_s \,,\, \nabla_t B_s \ra_{L^2(\Si)} \bigr)\bigr|_{t=0}\\
&= \la B_s \,,\, - \nabla_s B_t + \rd_A * F_A \ra_{L^2(\Si)}\bigr|_{t=0} \\
&= \int_\Si \la B_s \wedge 
                    \bigl( \nabla_s B_s - * (\pd_s *) B_s \bigr)\ra\bigr|_{t=0}  \\
&\leq \Bigl( C \bigl\| B_s \bigr\|_{L^2(\Si)}^2 
            - \int_\Si \la \nabla_s B_s \wedge B_s \ra \Bigr)\Bigr|_{t=0} .
\end{align*}
The second estimate for the normal derivative can be checked in any gauge at a fixed 
$(s_0,0)\in D\cap\pd\H^2$.
We choose a gauge with $\P\equiv 0$ and hence $B_s=\pd_s A$.
Then for the path $\X|_{(\cdot,0)\times\S}=A(\cdot,0)$ in $\cL_Y$ lemma~\ref{lem b}~(i) provides a
path of extensions $\tA : (s_0-\ep,s_0+\ep)\to\cA_{\rm flat}(Y)$ such that
$\pd_s\tA(s_0)|_\S = \pd_s A(s_0,0)$ and 
$\|\pd_s\tA(s_0)\|_{L^3(Y)} \leq C \|\pd_s A(s_0,0)\|_{L^2(\S)}$.
Here we fix a smooth path of metrics on $Y$ that extend the metrics $g_{s,0}$
on $\S$ for $s\in [-r_0,r_0]$.
The constant $C$ can then be chosen uniform for all $(s_0,0)\in D\cap\pd\H^2$.
So we calculate at $s=s_0$
\begin{align*}
 - \int_\Si \la \nabla_s B_s \wedge B_s \ra  
&=  \int_{\pd Y} \la \pd_s \tA \wedge \pd_s\pd_s \tA \ra  \\
&=  \int_Y \la \rd_\tA \pd_s \tA \wedge \pd_s^2 \tA \ra  
  - \int_Y \la \pd_s \tA \wedge \rd_\tA\pd_s^2 \tA \ra  \\
&=  \int_Y \la \pd_s \tA \wedge [\pd_s\tA \wedge \pd_s\tA] \ra  \\
&\leq \|\pd_s\tA\|_{L^3(Y)}^3 
\,\leq\, C^3 \|\pd_s A\|_{L^2(\S)}^3
\,=\, C^3 \|B_s\|_{L^2(\S)}^3 .
\end{align*}
Here we have used the fact that $F_\tA\equiv 0$, so $\rd_\tA\pd_s\tA = \pd_s F_\tA = 0$ 
and
$$
0 \;=\; \pd_s^2 F_\tA \;=\; \rd_\tA\pd_s^2\tA + [\pd_s\tA \wedge \pd_s\tA] .
$$

\vspace{-5.5mm}

\hspace*{\fill}$\Box$
\medskip

The significance of the following lemma is that a uniform bound on the slicewise 
$L^2$-norm of the curvature of an anti-self-dual connection implies an $L^p$-bound
on the curvature for any $p<3$. The specific value of the latter bound is not relevant
here. We only give it for comparison with a similar calculation in the proof of
proposition~\ref{prp est}.

\begin{lem} \label{lem L2 bound}
Fix $r_0>0$, let $2<p<3$, and let $\cm$ be a $\cC^1$-compact 
set of metrics of normal type on $D\times\Si$.
Then there exists a constant $C_p$ such that the following holds for all 
$0<\d\leq\half r_0$.

Let $\X\in\cA(D_{2\d}(0)\times\S)$ be anti-self-dual with respect to a metric in $\cm$
and suppose that for some constant $c$
$$
\bigl\| F_\X (s,t) \bigr\|_{L^2(\Si)} \leq c \qquad\forall (s,t)\in D_{2\d}(0).
$$
Then
$$
\bigl\| F_\X \bigr\|_{L^p(D_\d(0)\times\Si)} 
\;\leq\; C_p \bigl( \d^{\frac 4p -1} \, c + \d^{\frac 2p} \, c^{p-\frac 2p} \bigr).
$$
\end{lem}
\Pr
Fix a metric of normal type on $D\times\Si$. It suffices to prove the 
estimate with a uniform constant for all metrics of normal type in a 
$\cC^1$-neighbourhood of the fixed metric.
We choose this neighbourhood such that we have a uniform constant $C_1$ 
in the estimate from lemma~\ref{lem Laplacians},
$$
\laplace \bigl| F_\X \bigr|^2 
\leq  C_1 \bigl| F_\X \bigr|^2 + 8 \bigl| F_\X \bigr|^3 .
$$
Next, the normal coordinates at any $(s,t,z)\in D_{\frac 12 r_0}(0)\times\Si$ 
give a coordinate chart on $B_R(s,t,0,0)\cap\H^4$ with $R>0$ in which the fixed
metric (and hence all metrics in a sufficiently small neighbourhood) is 
$\cC^1$-close to the Euclidean metric. This $R>0$ can be chosen uniform for all
$(s,t,z)\in D_{\frac 12 r_0}(0)\times\Si$ such that the metrics in the coordinates 
meet the assumption of proposition~\ref{prp mean int}.
Now let $\bm:=\frac\m{64}$ where $\m>0$ is the constant from the theorem, and 
assume that $(s,t,z)\in D_\d(0)\times\Si$.
One can then apply this mean value inequality to $e=| F_\X |^2$ on 
$B_r(s,t,0,0)$ for $r=\min\bigl(t,R,c^{-1}\sqrt{\bm/\pi}\,\bigr)$.
Since
$$
\int_{B_r(s,t,0,0)} \bigl| F_\X \bigr|^2 
\;\leq\; \int_{B_r(s,t)} \bigl\| F_\X \bigr\|_{L^2(\Si)}^2 
\;\leq\; \pi \, r^2 c^2 \;\leq\; \bm
$$ 
we obtain with a uniform constant $C$ for all $(s,t,z)\in D_\d(0)\times\Si$
$$
\bigl| F_\X(s,t,z) \bigr|^2 
\;\leq\; C r^{-4} \int_{B_r(s,t,0,0)} \bigl| F_\X \bigr|^2 
\;\leq\; C \pi c^2 \min\bigl(t,R,c^{-1}\sqrt{\bm/\pi}\,\bigr)^{-2}.
$$
(Here we have used the fact that $r\leq R$, so 
${C_1}^2 + r^{-4} \leq C r^{-4}$ with a uniform constant depending on $R$.)
This pointwise control of $F_\X$ combines with the bound on
$\| F_\X (s,t) \|_{L^2(\Si)}$ to yield for $2<p<3$
\begin{align*}
\bigl\| F_\X \bigr\|_{L^p(D_\d(0)\times\Si)}^p
&\leq  \int_{D_\d(0)} \bigl\| F_\X \bigr\|_{L^\infty(\Si)}^{p-2}
                                \bigl\| F_\X \bigr\|_{L^2(\Si)}^2  \ds\dt  \\
&\leq \d c^2 
      \int_0^\d ( C\pi c^2 )^{\frac{p-2}2} 
        \bigl( t^{2-p} + \min\bigl(R,c^{-1}\sqrt{\bm/\pi}\,\bigr)^{2-p} \bigr) \;\dt \\
&\leq ( C\pi )^{\frac{p-2}2}  \d c^p \Bigl( \tfrac 1{3-p} \d^{3-p} 
                 + \d R^{2-p} + \d c^{p-2}\bigl(\tfrac\bm\pi\bigr)^{\frac{2-p}2} \Bigr) \\
&\leq {C_p}^p \bigl( \d^{4-p} \, c^p + \d^2 \, c^{2p-2} \bigr).
\end{align*}

\vspace{-5.5mm}

\hspace*{\fill}$\Box$
\medskip

Note that the assumption $p<3$ is crucial in this estimate. So a pointwise blowup
of the curvature is not enough to deduce a blowup of $\| F_\X \|_{L^2(\Si)}$.
As a consequence, it is essential that the compactness result 
\cite[Thm~B]{W elliptic} for solutions of (\ref{bvp}) with an $L^p$-bound on the
curvature was established for $2<p\leq 4$ (as well as for the easier case $p>4$).
These results put us in the following position near any slice of the boundary: 
There either is a local $L^p$-bound with $2<p<3$ for the curvature 
(and hence a convergent subsequence up to gauge) or a blowup of the functions
$\| F_{\X^\n} \|_{L^2(\Si)}^2: D\to [0,\infty)$.

If one now tries to mimic the proof of theorem~\ref{thm quant0}, one firstly needs
the following mean value inequality for the Laplacian with Neumann boundary condition, 
a proof of which can be found in \cite{W mean}.

\begin{prp}
\label{prp mean bdy}
There exist constants $C$, $\m>0$ such that the following holds.

Let $D_r(y)\subset\H^2$ be a Euclidean ball of radius $r>0$ and center 
$y\in\H^2$ intersected with the half space. 
Suppose that $e\in\cC^2(D_r(y),[0,\infty))$ satisfies for some constants
$A,B\geq 0$
\begin{equation*}
\left\{\begin{array}{ll}
\laplace e &\leq  B e , \\
-\tfrac\pd{\pd t}\bigr|_{\pd\H^2} e \hspace{-2mm} 
&\leq A \bigl( e + e^{\frac 32} \bigr) ,
\end{array}\right.
\qquad\text{and}\qquad
\int_{D_r(y)} e < \m A^{-2} .
\end{equation*}
Then
$$
e(y) \leq C \bigl( B + A^2 + r^{-2} \bigr) \int_{D_r(y)} e .
$$
\end{prp}

Another ingredient in our proof of energy quantization is the Hofer trick, 
\cite[6.4~Lemma~5]{HZ}, which we state here for convenience.

\begin{lem} {\bf (Hofer trick)} \label{lem HT}
Let $f:X\to[0,\infty)$ be continuous on the complete metric space $X$.
Then for every $x_0\in X$ and $\ep_0>0$ there exist $x\in B_{2\ep_0}(x_0)\subset X$
and $0<\ep\leq\ep_0$ such that $\ep f(x) \geq \ep_0 f(x_0)$ and
$f(y) \leq 2 f(x)$ for all $y\in B_\ep(x)$.
\end{lem}

The assumptions of proposition~\ref{prp mean bdy} will be verified by lemma~\ref{lem Laplacians}.
Firstly, the estimate for the normal derivative at the boundary, 
$-\tfrac\pd{\pd t}\bigr|_{t=0} e^\n$, results from lemma~\ref{lem b}~(i), i.e.\ from a (linear) 
extension of tangent vectors to $\cL_Y$ to \hbox{$1$-forms} on $Y$. 
Secondly, one should note that the term 
$\la F_A \,,\, [B_s\wedge B_s] \ra_{L^2(\Si)}$ in the expression for 
$\laplace \| F_\X \|_{L^2(\Si)}^2$ in lemma~\ref{lem Laplacians} 
is not yet in a form that can be controlled by any power of
$\| F_\X \|_{L^2(\Si)}^2$ as required above.
This is the central analytic problem of the bubbling analysis.
It is overcome by the following proposition which shows that 
$\| F_A \|_{L^\infty(\Si)}$ is essentially bounded by 
$\| F_\X \|_{L^2(\Si)}^2$.

If this bound was not true, then one would roughly find a pointwise blowup of the
$F_A$-component of the curvature while the energy goes to zero. 
A local rescaling would then lead to a nonflat limit connection in contradiction 
to the vanishing of the energy.
The nontrivial limit is obtained only when the blowup is mainly in the $F_A$-component
of the curvature. 
This is since after the local rescaling one has $\cC^0$-convergence only for $F_A$ 
(which satisfies a Dirichlet boundary condition) and not for $B_s$ 
(for which the global Lagrangian boundary condition is lost).
We will first state this result and show how it leads to a proof of theorem~\ref{thm A},
and then give its actual proof.

Recall that the boundary value problem (\ref{bvp}) is the anti-self-duality equation
together with a Lagrangian boundary condition in the space of flat connections over 
$\Si$. For the proposition below, it would actually suffice to assume only the flat
boundary condition $F_\X|_{(s,0)\times\Si} = 0$ in (\ref{bvp}).

\begin{prp} \label{prp est} \hspace{1mm} \\
\vspace{-5mm} 
\begin{enumerate}
\item 
Let $\X^\n\in\cA(D\times\Si)$ be a sequence of solutions of (\ref{bvp}) 
such that for some $\half r_0\geq\d>0$
$$
\sup_\n \sup_{(s,t)\in D_{2\d}} \bigl\| F_{\X^\n} (s,t) \bigr\|_{L^2(\Si)}  < \infty .
$$
Then
$$
\sup_\n \sup_{(s,t)\in D_\d} \bigl\| F_{A^\n} (s,t) \bigr\|_{L^\infty(\Si)} 
 < \infty .
$$
\item
For every $\cC^3$-compact set $\cm$ of metrics of normal type on $D\times\Si$
and every $\D>0$ there exists a constant $C$ such that the following holds:

Let $\X^\n\in\cA(D\times\Si)$ be a sequence of solutions of (\ref{bvp}) with
respect to metrics $g^\n\in\cm$. 
Suppose that $(s^\n,t^\n)\in D_{\frac 12 r_0}$, $\ep^\n\to 0$, and $R^\n\to\infty$ 
such that $\ep^\n R^\n\geq\D>0$ for all $\n\in\N$ and
$$
\quad\bigl\| F_{\X^\n} (s,t) \bigr\|_{L^2(\Si)}  \leq R^\n  
\qquad\qquad\forall (s,t)\in D_{2\ep^\n}(s^\n,t^\n) .
$$
Then for all sufficiently large $\n\in\N$
$$
\bigl\| F_{A^\n} (s,t) \bigr\|_{L^\infty(\Si)}   \leq  C \bigl(R^\n\bigr)^2 
\qquad\forall (s,t)\in D_{\ep^\n}(s^\n,t^\n).
$$
\end{enumerate}
\end{prp}

\noindent
{\bf Proof of theorem \ref{thm A}: } \\
Let $\cm$ be a $\cC^3$-compact set of metrics of normal type on $D\times\Si$ and
consider a sequence $\X^\n\in\cA(D\times\Si)$ of solutions of (\ref{bvp}) with respect
to metrics $g^\n\in\cm$.
We suppose that for some $2<p<3$ there is no local $L^p$-bound on the curvature near 
$\{0\}\times\Si$. 
By lemma~\ref{lem L2 bound} one then finds a subsequence (still denoted 
$(\X^\n)_{\n\in\N}$) and $D\ni (\bar s^\n,\bar t^\n)\to 0$ such that
$\bar R^\n :=\| F_{\X^\n}(\bar s^\n,\bar t^\n) \|_{L^2(\Si)} \to\infty$.
We pick $\bep^\n\to 0$ such that still $\bep^\n\bar R^\n\to\infty$.
The Hofer trick, lemma~\ref{lem HT}, then yields $D\ni(s^\n,t^\n)\to 0$ 
and $\ep^\n\to 0$ such that
$\| F_{\X^\n}(s^\n,t^\n) \|_{L^2(\Si)}=R^\n$ with $\ep^\n R^\n \to \infty$ and
$$
\bigl\| F_{\X^\n}(s,t) \bigr\|_{L^2(\Si)} \leq 2 R^\n 
\qquad\forall (s,t)\in D_{2\ep^\n}(s^\n,t^\n).
$$
Next, proposition~\ref{prp est}~(ii) asserts that for all $\n\geq\n_0$
$$
\bigl\| F_{A^\n} (s,t) \bigr\|_{L^\infty(\Si)} \leq  C \bigl(R^\n\bigr)^2 
\qquad\forall (s,t)\in D_{\ep^\n}(s^\n,t^\n).
$$
Here and in the following $C$ denotes any uniform constant.
Now consider the functions $e^\n= \| F_{\X^\n} \|^2 : D \to [0,\infty)$.
Use lemma~\ref{lem Laplacians} and the above bound to see that these satisfy 
on $D_{\ep^\n}(s^\n,t^\n)$ 
$$
\laplace e^\n
\;\leq\; C e^\n - 20 \la F_{A^\n} \,,\, [B^\n_s\wedge B^\n_s] \ra_{L^2(\Si)} 
\;\leq\; C \bigl( 1 + (R^\n)^2 \bigl) e^\n .
$$
For the normal derivative we obtain from lemma~\ref{lem Laplacians}
with a uniform constant $A$
$$
- \tfrac\pd{\pd t}\bigr|_{t=0} e^\n
\;\leq\; C e^\n - 4\int_\Si \la \nabla_s B^\n_s \wedge B^\n_s \ra  
\;\leq\; A \bigl( e^\n + (e^\n)^{\frac 32} \bigr) .
$$
Next, let $\m>0$ be the constant from the mean value inequality proposition~\ref{prp mean bdy}.
Now if $\n\geq\n_0$ and
\begin{equation}\label{e small}
\int_{D_{\ep^\n}(s^\n,t^\n)} e^\n \leq \m A^{-2}
\end{equation}
then this proposition asserts that
$$
(R^\n)^2 \;=\; e^\n(s^\n,t^\n) \;\leq\; 
C \bigl( (R^\n)^2 + (\ep^\n)^{-2} + 1 \bigr) \int_{D_{\ep^\n}(s^\n,t^\n)} e^\n .
$$
From this it would follow that
$$
\int_{D_{\ep^\n}(s^\n,t^\n)} e^\n \;\geq\; 
C^{-1} \bigl( 1 + (\ep^\n R^\n)^{-2} + (R^\n)^{-2} \bigr)^{-1} .
$$
Hence for all $\n\geq\n_0$ we must either have (\ref{e small}) violated or this 
inequality holds true. 
Now for sufficiently large $\n$ the right hand side is bounded below by 
$\half C^{-1}$, thus in any case
$$
\int_{D_{\ep^\n}(s^\n,t^\n)\times\Si} \bigl| F_{\X^\n} \bigr|^2 
\;>\; \min (\half C^{-1}, \m C^{-2} ) \;=:\, \ep_0 .
$$

\vspace{-5.5mm}

\hspace*{\fill}$\Box$
\medskip

The proof of proposition~\ref{prp est} is based on the following boundary
regularity result for anti-self-dual instantons on the half space with 
slicewise flat boundary conditions. These will arise from a local 
rescaling construction.

Here we use the coordinates $(x,y,s,t)$ with $t\geq 0$ on $\H^4$,
and as before we write connections $\X\in\cA(\H^4)$ in the splitting
$\X = A + \P\ds + \Psi\dt$.
Note that under the assumptions of the following lemma (with any $p>2$), 
the strong Uhlenbeck compactness for anti-self-dual connections (e.g.\ \cite[Thm~E]{W}) 
immediately implies the $\cC^\infty$-convergence of a subsequence of connections 
(in a suitable gauge) in the interior, away from $\pd\H^4$.
The slicewise flat boundary conditions are not quite enough to also obtain this convergence
at the boundary, however we still obtain some partial regularity results for this boundary
value problem. These provide the $\cC^0$-convergence of the curvature component $F_A$,
that vanishes on the boundary.

\begin{lem} \label{lem ASD on H}
Let $p>\frac 83$ and let $D_1(0)\subset\H^4$ be the unit half ball of radius $1$.
Let $g^\n$ be a sequence of metrics on $D_1(0)$ that converges to the Euclidean 
metric in the $\cC^3$-norm.
Let $\X^\n\in\cA^{1,p}(D_1(0))$ be a sequence of anti-self-dual connections with respect
to the metrics $g^\n$ and that satisfy flat boundary conditions, $F_{A^\n}|_{t=0}\equiv 0$.
Suppose that 
$$
\lim_{\n\to\infty} \bigl\| F_{\X^\n} \bigr\|_{L^p(D_1(0))} = 0 .
$$
Then there exists a subsequence such that
$$
\lim_{\n\to\infty} \bigl\| F_{A^\n} \bigr\|_{L^\infty(D_{\frac 12}(0))} = 0 .
$$
\end{lem}
\Pr 
Let $U\subset\H^4$ be a compact submanifold with smooth boundary obtained from
$D_{\frac 34}(0)$ by 'rounding off the edge' at $\pd\H^4$, so 
$D_{\frac 12}(0)\subset U \subset D_1(0)$.
(More precisely, Uhlenbeck's gauge theorem below requires that the domain $U$ is diffeomorphic
to a ball; to obtain uniform constants, it should moreover be starlike w.r.t.\ $0$.)
Let a sequence of connections $\X^\n$ as above be given.
For sufficiently large $\n$ the metrics $g^\n$ on $U$ are all sufficiently $\cC^2$-close to 
the Euclidean metric so that the Uhlenbeck gauge for $\X^\n|_U$ exists with uniform
constants: The energy $\int_U |F_{\X^\n}|^2$ becomes arbitrarily small for large $\n$, so by
\cite[Thm~1.3]{U2} or \cite[Thm~6.3]{W} these connections can be put into a gauge 
(again denoted by $\X^\n\in\cA^{1,p}(U)$) such that
$\rd^*\X^\n =0$ and $*\X^\n|_{\pd U}=0$.
%
%
%
%
This gauge also gives a uniform bound 
$\|\X^\n\|_{W^{1,p}(U)}\leq C_{Uh}\|F_{\X^\n}\|_{L^p(U)}\leq C$.

We now have to follow through the higher regularity arguments in the proof of
\cite[Thm~2.6]{W elliptic} to find a uniform bound on $F_{A^\n}$ in the
H\"older norm $\cC^{0,\l}(D_{\frac 12}(0))$ for some $\l>0$.
This will finish the proof since the embedding $\cC^{0,\l}\hookrightarrow\cC^0$
is compact, so this would imply $\cC^0$-convergence for a subsequence.
The limit can only be $0$ since that was already the $L^p$-limit on $D_1(0)$.

Firstly note that the metrics on $U$ are all $\cC^3$-close, so 
\cite[Thm~2.6]{W elliptic} allows for uniform estimates up to the third 
derivatives of the connections.\footnote
{
The original theorem requires $\cC^5$-bounds, but $\cC^3$-bounds suffice when the metrics are 
already given in the appropriate coordinates (that otherwise are determined by the metrics).
}
Since $2<p<3$ we are dealing with 'Boundary case II' in the proof of this theorem.

We choose $d>\half$ such that $Q_d:=[-d,d]\times[0,d]\times B_d\subset U$, where 
$B_d\subset\R^2$ is the Euclidean ball centered at $0$.
We moreover drop the superscript $\n$. Then the connections 
$\X=A+\P\ds+\Psi\dt$ with $A:[-d,d]\times[0,d]\to\Om^1(B_d;\cg)$ and
$\P,\Psi:[-d,d]\times[0,d]\to\Om^0(B_d;\cg)$ solve the following
boundary value problem analogous to \cite[(12)]{W elliptic}.\footnote
{
Note that $B$ is replaced by $A$ and we use the reference connection $A_0=0$.
}
Here $Q_d$ is equipped with the metric $\ds^2+\dt^2 +g_{s,t}$, and we shall write
$\rd$, $\rd^*$ and $\nabla$ for the families of operators on $B_d$ with respect 
to the metrics $g_{s,t}$.
Note that due to the localization we only retain the flat boundary condition.
\begin{equation}\label{bvp comp}
\left\{\begin{aligned} 
\rd^*A &= \pd_s\P + \pd_t\Psi , \\
*F_A &= \pd_t\P - \pd_s\Psi  + [\Psi,\P] ,\\
\pd_s A + * \pd_t A &= \rd_A\P + * \rd_A\Psi , \\
\Psi(s,0) =0 &\quad\forall s\in[-d,d] , \\
F_A(s,0) =0 &\quad\forall s\in[-d,d] .
\end{aligned}\right.
\end{equation}
Firstly, this combines to Laplace equations on $\P$ and $\Psi$ (see (\ref{Phi lap}) 
below) with a Dirichlet boundary condition for $\Psi$ and an inhomogeneous Neumann 
condition for $\P$,
$$
\pd_t\P|_{t=0} = \pd_s\Psi - [\Psi,\P] .
$$
By e.g.\ \cite[Prop~2.7]{W elliptic} this yields $W^{2,q}$-bounds for
$\P$ and $\Psi$ on $Q_d$ with a slightly smaller $d>\half$. Due to nonlinearities in
the lower order terms, these bounds hold only for $q=\frac{4p}{8-p}$ (i.e.\ when
$W^{1,p}\cdot L^p\hookrightarrow L^q$). However, we have assumed $p>\frac 83$ so that
$q>2$ and thus $W^{2,q}(Q_d)$ embeds into $\cC^0(Q_d)$.

Next, one has $W^{1,q}$-bounds on $\rd^*A$ and $\rd A$ as in \cite[(13)]{W elliptic}.
These lead to a bound on $\nabla A \in W^{1,q}(Q_d)$ (again for slightly smaller 
$d>\half$), see \cite[Lemma~2.9]{W elliptic}.
In particular, $A$ is bounded in
$W^{1,q}([-d,d]\times[0,d],W^{1,q}(B_d))$, which embeds into $\cC^0(Q_d)$.
Thus we have obtained $\cC^0$-bounds on the whole connection $\X$.
Using these in the nonlinear terms and going through the previous two steps again
yields bounds on $\P,\Psi\in W^{2,p}(Q_d)$ and $\nabla A\in W^{1,p}(Q_d)$ (with
slightly smaller $d>\half$).
In order to obtain bounds for third derivatives of $\P$ and $\Psi$ we calculate 
\begin{align}
\laplace\P 
&= \pd_s \bigl( \pd_t\Psi - \rd^*A \bigr) 
  +\pd_t \bigl( [\P,\Psi] - \pd_s\Psi - *F_A \bigr) \nonumber \\
&\quad  +\rd^* \bigl( \pd_s A + *\pd_t A + *\rd_A\Psi - [A,\P] \bigr) \label{Phi lap}  \\
&= \pd_t \bigl( [\P,\Psi] - *[A\wedge A] \bigr) 
   -\rd^* [A,\P] - *\rd[A,\Psi] + \text{l.o.} \;. \nonumber 
\end{align}
Here we have disregarded all lower order terms that arise from derivatives of the 
metric. From this one obtains an $L^q$-bound on $\laplace\nabla\P$. Indeed, the 
crucial terms are $*[A \wedge \nabla\pd_t A]$ and $*[\nabla A \wedge \pd_t A]$,
where in both cases the first factor is $W^{1,p}$-bounded and the second factor
is $L^p$-bounded. The analogous calculation also works for $\Psi$, so with the boundary
conditions as before we obtain (for smaller $d>0$) bounds on
$\nabla\P,\nabla\Psi\in W^{2,q}(Q_d)$. 
In particular this gives bounds for $\P$ and $\Psi$ in 
$W^{2,q}([-d,d]\times[0,d],W^{1,q}(B_d))$, and thus
$*F_A=\pd_t\P-\pd_s\Psi+[\Psi,\P]$ is bounded in
$W^{1,q}([-d,d]\times[0,d],W^{1,q}(B_d))$.
Now finally, there is a continuous embedding of $W^{1,q}$ (on a $2$-dimensional domain
with values in any Banach space) into the H\"older space $\cC^{0,2\l}$ with some
$\l>0$, so the above space embeds into 
$\cC^{0,2\l}([-d,d]\times[0,d],\cC^{0,2\l}(B_d))$, 
which in turn embeds continuously into $\cC^{0,\l}(Q_d)$.
Thus we obtain the claimed uniform bounds on $F_{A^\n}\in\cC^{0,\l}(Q_d)$.
\QED \\

\noindent
{\bf Proof of proposition \ref{prp est}: }
 (i) and (ii) are proven by the same contradiction.

If (ii) was not true, then one would have a $\cC^3$-compact set $\cm$ of metrics of 
normal type on $D\times\Si$ and $\D>0$ with the following significance.
For all $k\in\N$ there is a sequence 
$\X_k^\n\in\cA(D\times\Si)$ of solutions of (\ref{bvp}) with respect to metrics 
$g_k^\n\in\cm$, moreover $(\bar s_k^\n,\bar t_k^\n)\in D_{\frac 12 r_0}$, 
$\ep_k^\n\to 0$, and $R_k^\n\to\infty$ such that $\ep_k^\n R_k^\n\geq\D>0$ and

\vspace{-5mm}

$$
\qquad\qquad\qquad
\bigl\| F_{\X_k^\n} (s,t) \bigr\|_{L^2(\Si)}  \leq R_k^\n 
\qquad\forall(s,t)\in D_{2\ep_k^\n}(\bar s_k^\n,\bar t_k^\n) . 
$$
But for every $k\in\N$ and $\n_0\in\N$ there would exist $\n\geq\n_0$,
$(s_k,t_k)\in D_{\ep_k^\n}(\bar s_k^\n,\bar t_k^\n)$, and $z_k\in\Si$
such that
$$
\bigl| F_{A_k^\n} (s_k,t_k,z_k) \bigr| >  k \bigl(R_k^\n\bigr)^2 .
$$
For each $k\in\N$ we choose $\n_0$ such that $\ep^\n_k\leq \frac 1k$ and $R^\n_k\geq k$
for all $\n\geq\n_0$. 
Then from a subsequence of a diagonal sequence one obtains
solutions $\X^k\in\cA(D\times\S)$ of (\ref{bvp}) w.r.t.\ a $\cC^3$-convergent
sequence of metrics $g^k\to g^\infty$; 
$\ep^k\to 0$ and $R^k\to\infty$ such that $\ep^k R^k \geq \D$;
$(s^k,t^k)\to (s^\infty,t^\infty)\in D_{\frac 12 r_0}$, 
$\Si\ni z^k \to z$, and $C^k\to\infty$ such that

\vspace{-5mm}

$$
\bigl| F_{A^k} (s^k,t^k,z^k) \bigr| \geq  \bigl(C^k R^k\bigr)^2 .
$$
Since $D_{\ep_k^\n}(s_k,t_k)\subset  D_{2\ep_k^\n}(\bar s_k^\n,\bar t_k^\n)$
one also obtains the bound
$$
\bigl\| F_{\X^k} (s,t) \bigr\|_{L^2(\Si)} \leq R^k 
\qquad\forall (s,t)\in D_{\ep^k}(s^k,t^k).
$$
If (i) was not true, then one would find a sequence $\X^\n\in\cA(D\times\Si)$ 
solving (\ref{bvp}), constants $C$, $0<\d\leq\half r_0$, and 
$(s^\n,t^\n,z^\n)\in D_\d(0)\times\Si$ such that 
$|F_{A^\n}(s^\n,t^\n,z^\n)|\to\infty $ but for all $\n\in\N$
$$
\sup_{(s,t)\in D_{2\d}(0)} \bigl\| F_{\X^\n} (s,t) \bigr\|_{L^2(\Si)} \leq C .
$$
For a subsequence we can assume that 
$(s^\n,t^\n,z^\n)\to (s^\infty,t^\infty,z^\infty)\in D_{\frac 12 r_0}\times\S$.
We set $R^\n=C^\n:=|F_{A^\n}(s^\n,t^\n,z^\n)|^{\frac 14}\to\infty$
and $\ep^\n:=\min\bigl( (R^\n)^{-1}, \d \bigr)\to 0$. 
Then one has $C\leq R^\n$, $\ep^\n R^\n \geq 1=:\D$, and 
$D_{\ep^\n}(s^\n,t^\n)\subset D_{2\d}(0)$ for all sufficiently large $\n\in\N$.
That way, we have constructed the same sequences as for (ii), to which we shall
find a contradiction:
\begin{itemize}
\item 
Solutions $\X^\n\in\cA(D\times\Si)$ of (\ref{bvp}) with respect to 
$\cC^3$-convergent metrics $g^\n\to g^\infty$, constants
$\ep^\n\to 0$, $R^\n\to\infty$ with $\ep^\n R^\n \geq\D > 0$, and $C^\n\to\infty$,
and points $(s^\n,t^\n,z^\n)\to(s^\infty,t^\infty,z^\infty)\in D_{\frac 12 r_0}\times\S$ such that
$$
\sup_{(s,t)\in D_{\ep^\n}(s^\n,t^\n)} \bigl\| F_{\X^\n} (s,t) \bigr\|_{L^2(\Si)} 
\leq R^\n , \quad\qquad
\bigl| F_{A^\n} (s^\n,t^\n,z^\n) \bigr| \geq  \bigl(C^\n R^\n\bigr)^2 . 
$$
\end{itemize}
Firstly suppose that $\lim\sup_\n t^\n R^\n \geq d>0 $.
In that case choose $d>0$ even smaller so $\d\leq \D$, then $0<r^\n:=d(R^\n)^{-1}\leq\ep^\n$ and
$r^\n\leq t^\n$ for a suitable subsequence.
Now the geodesic ball $B_{r^\n}(s^\n,t^\n,z^\n)$ with respect to $g^\n$ is entirely contained in 
$D\times\Si$, and for sufficiently large $\n$ it will be small enough to lie 
within a normal coordinate chart around $(s^\infty,t^\infty,z^\infty)$ for the metric 
$g^\infty$.
In this coordinate chart all metrics $g^\n$ for large $\n$ will be sufficiently 
$\cC^1$-close to the Euclidean metric for proposition~\ref{prp mean int} to apply with
uniform constants $\m>0$ and $C$.
Next, lemma~\ref{lem Laplacians} gives a uniform constant $C_1$ such that
\begin{equation}  \label{C1 Laplace}
\laplace \bigl| F_{\X^\n} \bigr|^2
\;\leq\; C_1 \bigl| F_{\X^\n} \bigr|^2 + 8 \bigl| F_{\X^\n} \bigr|^3 . 
\end{equation}
Let $\bm:=\frac\m{64}$ and choose $d\leq\sqrt{\bm/\pi}$ then
$$
\int_{B_{r^\n}(s^\n,t^\n,z^\n)} \bigl| F_{\X^\n} \bigr|^2 
\;\leq\; \pi (r^\n)^2 \bigl(R^\n\bigr)^2 \;\leq\; \bm .
$$
So proposition~\ref{prp mean int} implies
$$
\bigl(C^\n R^\n\bigr)^4 
\;\leq\; \bigl| F_{\X^\n} (s^\n,t^\n,z^\n) \bigr|^2
\;\leq\; C \bigl( {C_1}^2 + (r^\n)^{-4} \bigr) 
           \int_{B_{r^\n}(s^\n,t^\n,z^\n)} \bigl| F_{\X^\n} \bigr|^2 .
$$
Putting in above estimate of the energy and $r^\n R^\n = d>0$ then leads to
the contradiction 
$$
\bigl(C^\n \bigr)^4
\;\leq\; C \bm \bigl( {C_1}^2 \bigl(R^\n\bigr)^{-4}  + d^{-4} \bigr) 
\;\underset{\n\to\infty}{\longrightarrow}\; C \bm d^{-4} \;<\; \infty.
$$
The second and crucial case is $t^\n R^\n \to 0$.
We choose $\D\geq d>0$ and set $r^\n := d(R^\n)^{-1} \leq \ep^\n$ such that 
$\frac 13 r^\n \geq t^\n$ for sufficiently large $\n$.
Now for all $(s,t)\in D_{\frac 13 r^\n}(s^\n,t^\n)$ we have 
$t\leq t^\n + \frac 13 r^\n \leq \frac 23 r^\n$, hence 
$B_t(s,t)\subset D_{\ep^\n}(s^\n,t^\n)$, and thus for all $z\in\Si$

\pagebreak
$$
\int_{B_t(s,t,z)} \bigl| F_{\X^\n} \bigr|^2 
\;\leq\; \pi t^2 \bigl(R^\n\bigr)^2
\;\leq\;  \tfrac 49 \pi d^2 .
$$
As in the first case one can choose $\n$ sufficiently large such that for all $z\in\Si$
the above balls $B_t(s,t,z)\subset D_{\ep^\n}(s^\n,t^\n)\times\Si$ lie within
a normal coordinate chart around $(s^\infty,t^\infty,z)$ for the metric $g^\infty$.
Again, for large $\n$ all metrics $g^\n$ in these coordinates will be sufficiently 
$\cC^1$-close to the Euclidean metric, so that (\ref{C1 Laplace}) holds with a 
uniform constant $C_1$ and proposition~\ref{prp mean int} applies with uniform constants 
$\m>0$ and $C$.
We choose $d>0$ sufficiently small so that $\tfrac 49 \pi d^2 \leq \tfrac\m{64}$,
then proposition~\ref{prp mean int} implies that for all 
$(s,t,z)\in D_{\frac 13 r^\n}(s^\n,t^\n)\times\Si$
$$
\bigl| F_{\X^\n}(s,t,z) \bigr|^2
\;\leq\; C \bigl( {C_1}^2 +  t^{-4} \bigr) \int_{B_t(s,t,z)} \bigl| F_{\X^\n} \bigr|^2
\;\leq\; C \pi \bigl( 1 +  t^{-2} \bigr) \bigl(R^\n\bigr)^2 .
$$
Note here that $C_1 t \leq C_1 (r^\n + t^\n) \leq 1$ for sufficiently large $\n$.
With the above pointwise control of the curvature we can interpolate similar to 
lemma~\ref{lem L2 bound} to find for any fixed $2<p<3$ and for all $r\leq\frac 13 r^\n$
\begin{align*}
\int_{D_r(s^\n,t^\n)\times\Si} \bigl| F_{\X^\n} \bigr|^p
&\leq \int_{D_r(s^\n,t^\n)} \bigl\| F_{\X^\n} \bigr\|_{L^\infty(\Si)}^{p-2}
                            \bigl\| F_{\X^\n} \bigr\|_{L^2(\Si)}^2  \\
&\leq C \bigl( R^\n \bigr)^p
\int_{D_r(s^\n,t^\n)} \bigl( 1 + t^{2-p} \bigr) \\
&\leq C \bigl( R^\n \bigr)^p \bigl( \pi r^2 + \tfrac {2r}{3-p}(t^\n+r)^{3-p} \bigr) 
\;\leq\; C \bigl( R^\n \bigr)^p (t^\n+r)^{4-p}
\end{align*}
Here $C$ denotes any uniform constant (depending on the choice of $p$).
Next, recall that $|F_{A^\n}(s^\n,t^\n,z^\n)|\geq(C^\n R^\n)^2$ and 
$\ep^\n C^\n R^\n \geq \D C^\n \to \infty$.
So by the usual local rescaling we can define connections $\tX^\n$ on increasingly 
large $4$-balls (intersected with half spaces) 
$B_{\ep^\n C^\n R^\n}(0) \cap \{t\geq -t^\n C^\n R^\n \}\subset\R^4$.
We use normal coordinates for $g^\infty$ near $(s^\infty,t^\infty,z^\infty)$ 
and write \hbox{$\R^4=\{(s,t,z)\st s,t\in\R, z\in\R^2\}$}, 
then these connections are defined by
$$
\tX^\n(s,t,z):=\X^\n\bigl((s^\n,t^\n,z^\n)+\tfrac 1{C^\n R^\n}(s,t,z) \bigr) .
$$
They satisfy the boundary condition $F_{\tA^\n}|_{t=-t^\n C^\n R^\n}=0$, 
and they are anti-self-dual with respect to the metrics
$\tg^\n(s,t,z):=g^\n((s^\n,t^\n,z^\n)+\tfrac 1{C^\n R^\n}(s,t,z))$.
Note that the coordinates were chosen such that on bounded domains the metric
$\tg^\infty$ (rescaled by $C^\n R^\n$) converges to the Euclidean metric in any $\cC^k$-norm.
Thus for large $\n$ the metrics $\tg^\n$ become arbitrarily $\cC^3$-close to 
the Euclidean metric.

Moreover, this construction is such that $|F_{\tA^\n}(0)|\geq 1$ for all $\n$.
On the other hand for all $\r>0$ we have $\r(C^\n R^\n)^{-1}\leq\frac 13 r^\n$ 
for sufficiently large $\n$ and thus
\begin{align}
\bigl\| F_{\tX^\n} \bigr\|_{L^p(B_\r(0)\cap \{t\geq -t^\n C^\n R^\n \})}^p
&=\; \bigl(C^\n R^\n\bigr)^{4-2p} 
   \int_{D_{\r(C^\n R^\n)^{-1}}(s^\n,t^\n,z^\n)} \bigl| F_{\X^\n} \bigr|^p \nonumber\\
&\;\leq C \bigl(C^\n R^\n\bigr)^{4-2p} \bigl(R^\n\bigr)^p (t^\n+\r(C^\n R^\n)^{-1})^{4-p} 
\nonumber\\
&\;\leq C \bigl(C^\n\bigr)^{4-2p} (t^\n R^\n + \r(C^\n)^{-1})^{4-p} \label{Lp to 0}\\
&\underset{\n\to\infty}{\longrightarrow}\; 0 .  \nonumber
\end{align}
If $\lim\sup_\n t^\n C^\n R^\n >0$, then we can choose a subsequence and $\r>0$ such that
$F_{\tX^\n}$ is defined on $B_\r(0)$ for all $\n$. Then the above estimate shows that
$|F_{\tX^\n}|\to 0$ in the $L^p$-norm on $B_\r(0)$
Due to the strong Uhlenbeck compactness for anti-self-dual connections 
(see e.g.\ \cite[Thm~E]{W}) one can find a subsequence and gauge transformations 
(which do not affect the norm of the curvature) 
such that this convergence is actually in the $\cC^0$-topology. 
This contradicts $|F_{\tA^\n}(0)|\geq 1$.

If $\t^\n:=t^\n C^\n R^\n \to 0$ then we need lemma~\ref{lem ASD on H} to obtain 
this contradiction.
We shift the connections $\tX^\n$ and metrics $\tg^\n$ by $\t^\n$ in the $t$-direction 
so they are defined on $D_{\ep^\n C^\n R^\n}(0,\t^\n,0)$. In particular, for sufficiently
large $\n$, they are all defined on $D_1(0)$.
Now the $\tX^\n$ satisfy flat boundary conditions at $t=0$, and they are anti-self-dual 
with respect to the shifted metrics $\tg^\n$.
Since the shifts $\t^\n$ converge to $0$, we moreover preserve the $\cC^3$-convergence of
the metrics $\tg^\n$ to the Euclidean metric.
By this shift we have $|F_{\tA^\n}(0,t^\n,0)|\geq 1$. 
Moreover, choose any $\frac 83<p<3$, then we have 
$\bigl\| F_{\tX^\n} \bigr\|_{L^p(D_1(0))}\to 0$ 
since for $\n$ sufficiently large $D_1(0)\subset D_2(0,\t^\n,0)$, 
and the latter is the domain in (\ref{Lp to 0}) after the shifting.
Now lemma~\ref{lem ASD on H} asserts that in fact $F_{\tA^\n}$ converges to $0$ in the $\cC^0$-norm 
on $D_{\frac 12}(0)$. 
This however contradicts the fact that $|F_{\tA^\n}(z^\n)|\geq 1$ for $z^\n=(0,t^\n,0)\to 0$.
\QED

\section{Extension of connections in $\mathbf{\cL_Y}$}
\label{extension}

This section is devoted to the proof of lemma~\ref{lem b}, that is to extension
constructions that relate connections in the Lagrangian $\cL_Y$ on $\S=\pd Y$ to 
flat connections on $Y$.
Throughout $Y$ is a handle body with boundary $\pd Y=\S$ a Riemann surface of genus $g$.
We moreover fix some $p>2$.
The Lagrangian $\cL_Y\subset\cA^{0,p}(\S)$ as introduced in \cite[Lemma~4.6]{W Cauchy} 
is given by
\begin{align*}
\cL_Y 
&= {\rm cl}_{L^p}\, \bigl\{ A\in\cA_{\rm flat}(\S) \st \exists 
                              \tA\in\cA_{\rm flat}(Y) : \tA|_\S=A \bigr\} \\
&= \bigl\{ u^*(A|_\Si) \st A\in\cA_{\rm flat}(Y), u\in\cG^{1,p}(\Si) \bigr\} \\
&= \bigl\{ A\in\cA^{0,p}_{\rm flat}(\Si) \st 
    \r_z(A) \in {\rm Hom}(\pi_1(Y,z),\cg) \subset {\rm Hom}(\pi_1(\Si,z),\rG)  \bigr\} .
\end{align*}
The space $\cA^{0,p}_{\rm flat}(\Si)$ of weakly flat $L^p$-connections was introduced in
\cite[Sec.3]{W Cauchy}.
If we fix any $z\in\S$, then every weakly flat connection is gauge equivalent to a smooth
connection via a gauge transformation in the based gauge group
$$
\cG^{1,p}_z(\S) = \bigl\{ u \in \cG^{1,p}(\S) \st u(z)=\one \bigr\} .
$$
Thus the based holonomy $\r_z$ is welldefined on $\cA^{0,p}_{\rm flat}(\S)$ by
first going to a smooth gauge and then calculating the holonomy along fixed generators
of $\pi_1(\S,z)$.
We moreover recall from \cite{W Cauchy} that $\cL_Y$ is a Banach submanifold of 
$\cA^{0,p}(\S)$, and it is Lagrangian with respect
to the symplectic form $\o(\a,\b)=\int_\S \la \a \wedge \b \ra$ in the sense that 
$\rT_A\cL_Y\subset L^p(\S,\rT^*\S\otimes\cg)$ is a maximal isotropic subspace for all
$A\in\cL_Y$. Finally, $\cL_Y$ has the structure of a $\cG^{1,p}_z(\Si)$-bundle 
over the $g$-fold product 
$\rG\times \cdots \times \rG \cong {\rm Hom}(\pi_1(Y,z),\rG)$,
$$
\cG^{1,p}_z(\Si) 
\;\hookrightarrow\;  \cL_Y
\;\overset{\r_z}{\longrightarrow}\; \rG^g .
$$
We will fix a bundle atlas by specifying local sections over a finite cover of $\rG^g$.
For that purpose we choose loops $\a_1,\b_1,\ldots,\a_g,\b_g \subset \Si$ disjoint from $z$ 
that represent the standard generators of $\pi_1(\Si)$ such that $\a_1,\ldots,\a_g$ generate 
$\pi_1(Y)$ and the only nonzero intersections are $\a_i \cap \b_i$.
One can then modify the $\a_i$ such that they run through $z$ and coincide in a neighbourhood
of $z$ but still do not intersect the $\b_j$ for $j\neq i$.

The based holonomy $\r_z :\cL_Y \to \rG^g \cong {\rm Hom}(\pi_1(Y,z),\rG)$ is now
given by the $g$ holonomies $\hol_{\a_i}:\cL_Y \to \rG$ for the paths $\a_i$ starting
and ending at $z$.

Next, we choose spanning discs of the $\b_i$ that are pairwise disjoint and intersect $\pd Y$
in $\b_i$ only. 
Their tubular neighbourhoods provide orientation preserving diffeomorphisms 
$\psi_i:[0,1]\times D \to Z_i$ (with $D\subset\R^2$ the unit disc) to disjoint 
neighbourhoods $Z_i\subset Y$ of the spanning discs.
They can be chosen such that $\a_i\cap Z_j = \emptyset$ for $i\neq j$ and such that
$\psi_i:[0,1]\times \{y\} \overset{\sim}{\to} \a_i\cap Z_i$ for some $y\in\pd D$.
We then fix the induced orientation for the $\a_i$.

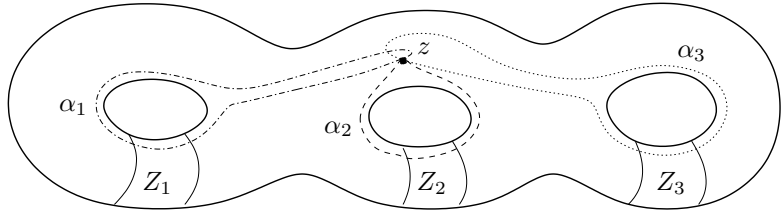
\begin{center}
\begin{picture}(0,0)%
\includegraphics{sigma.pstex}%
\end{picture}%
\setlength{\unitlength}{2368sp}%
\begingroup\makeatletter\ifx\SetFigFont\undefined%
\gdef\SetFigFont#1#2#3#4#5{%
  \reset@font\fontsize{#1}{#2pt}%
  \fontfamily{#3}\fontseries{#4}\fontshape{#5}%
  \selectfont}%
\fi\endgroup%
\begin{picture}(8122,2198)(1696,-2685)
\put(8724,-1069){\makebox(0,0)[lb]{\smash{\SetFigFont{10}{12.0}{\familydefault}{\mddefault}{\updefault}$\alpha_3$}}}
\put(8502,-2468){\makebox(0,0)[lb]{\smash{\SetFigFont{10}{12.0}{\familydefault}{\mddefault}{\updefault}$Z_3$}}}
\put(2251,-1636){\makebox(0,0)[lb]{\smash{\SetFigFont{10}{12.0}{\familydefault}{\mddefault}{\updefault}$\alpha_1$}}}
\put(5026,-1861){\makebox(0,0)[lb]{\smash{\SetFigFont{10}{12.0}{\familydefault}{\mddefault}{\updefault}$\alpha_2$}}}
\put(3124,-2468){\makebox(0,0)[lb]{\smash{\SetFigFont{10}{12.0}{\familydefault}{\mddefault}{\updefault}$Z_1$}}}
\put(5998,-2468){\makebox(0,0)[lb]{\smash{\SetFigFont{10}{12.0}{\familydefault}{\mddefault}{\updefault}$Z_2$}}}
\put(6001,-1036){\makebox(0,0)[lb]{\smash{\SetFigFont{10}{12.0}{\familydefault}{\mddefault}{\updefault}$z$}}}
\end{picture}

\end{center}

Choose $\D>0$ less than the injectivity radius of $\exp:\cg\to\rG$ and
fix a function $\t\in\cC^\infty([0,1],[0,1])$ with $\t|_{[0,\frac 14]}\equiv 0$ and
$\t|_{[\frac 34,1]}\equiv 1$.
Now given any fixed $\Th^0=(\th^0_1,\ldots,\th^0_g)\in\rG^g$ we choose smooth paths
$\g^0_i:[0,1]\to\rG$ with $\g^0_i|_{[0,\frac 14]}\equiv\one$ and 
$\g^0_i|_{[\frac 34,1]}\equiv\th_i^{-1}$.
Let $B_\D(\Th^0)\subset\rG^g$ be the closed exponential ball around $\Th^0$.
Then for every $\Th=(\th_i)\in B_\D(\Th^0)$ we have local gauge transformations
$v_i\in\cG(Z_i)$ given by $v_i(\psi_i(t,x))=( \g^0_i(t) \exp(\t(t)\x_i) )^{-1}$, where 
$\x_i=\exp^{-1}((\th^0_i)^{-1}\th_i)\in  B_\D(0)\subset\cg$.
Note that $v_i\equiv\one$ and $v_i\equiv\th_i^{-1}$ near the two boundary components of $Z_i \subset Y$.
These local gauge transformations can be used to define a local section of $\cL_Y$, that is
a smooth map $\X: B_\D(\Th^0) \to \cA_{\rm flat}(Y)$ such that $\r_z(\X(\Th)|_{\S})=\Th$,
\begin{equation} \label{section}
\X(\Th) \;=\; \left\{ 
\begin{array}{ll}
v_i^{-1}\rd v_i &; \,\text{on}\;Z_i , \\
0 &; \,\text{on}\;Y \setminus \bigcup_{i=1}^g Z_i .
\end{array}
\right.
\end{equation}
We now fix $\Th^0_j\in\rG^g$ for $j=1,\ldots,N$ such that the domains
$ B_{\frac 12\D}(\Th^0_j)$ already cover all of $\rG^g$.
This gives rise to a bundle atlas for $\cL_Y$ given by the charts
\begin{equation} \label{chart}
\begin{array}{ccc}
\cG^{1,p}_z(\Si) \times  B_\D(\Th^0_j) &\longrightarrow& \cL_Y \\
(u,\Th) &\longmapsto& u^*\X_j(\Th)|_\S .
\end{array}
\end{equation}
Next, we can find tubular neighbourhoods 
$\ta_i:[-1,1]\times [0,1] \hookrightarrow \S$ of the loops \hbox{$\a_i=\ta_i(0,\cdot)$}
that again coincide near $z$ for all $i=1,\ldots,g$.
Then these are a family of loops based at $\ta_i(\t,0)=\ta_i(\t,1)=z(\t)$ for some $i$-independent
smooth path $z:[-1,1]\to\S\setminus\bigcup Z_i$. 
As before, the intersection $\ta_i(\t,\cdot)\cap Z_j$ will be empty for $i\neq j$, 
and for $i=j$ it is $\psi_i([0,1]\times\{y(\t)\}$ for some $y(\t)\in\pd D$.

Note that for the special connections $\X(\Th)\in\cA_{\rm flat}(Y)$ as in (\ref{section})
the holono\-mies $\hol_{\ta_i(\t)}(\X(\Th)) = \hol_{\a_i}(\X(\Th))$ are independent 
of $\t\in[-1,1]$.
For other connections, the variation of the paths $\ta_i(\t,\cdot)$ along $\t\in[-1,1]$
allows to control the holonomy by the connections in the $L^1$-topology.

\begin{lem}  \label{lem holonomy}
There exists a constant $C$ such that the following holds.
\begin{enumerate}
\item 
For all smooth paths $A:(-\ep,\ep)\to\cA(\S)$ 
there exists $\t\in[-1,1]$ such that 
with $\th=\hol_{\ta_i(\t)}(A(0))$ 
for all $i=1,\ldots,g$
$$
\bigl| \pd_s|_{s=0} \hol_{\ta_i(\t)}(A(s)) \bigr|_{\rT_\th\rG}
\;\leq\; C \| \pd_s A (0) \|_{L^1(\S)} .
$$
\item
For all $A_0,A_1\in\cA(\S)$ there exists $\t\in[-1,1]$ such that for all $i=1,\ldots,g$
$$
{\rm dist}_\rG\bigl( \hol_{\ta_i(\t)}(A_0) \,,\,  \hol_{\ta_i(\t)}(A_1)  \bigr)
\;\leq\; C \| A_0 - A_1 \|_{L^1(\S)} .
$$
\end{enumerate}
\end{lem}
\Pr
Starting with the proof of (ii) we recall that for every $i=1,\ldots,g$ and all $\t\in[-1,1]$ 
the holonomies $\hol_{\ta_i(\t)}(A_j)=u_j(1)\in\rG$ for $j=0,1$
are given by the solutions $u_j:[0,1]\to\rG$ of 
\[
\dot u_j u_j ^{-1} = - \ta_i(\t)^*A_j 
\qquad\text{with}\qquad u_j(0) = \one.
\]
Note that for fixed $i=1,\ldots,g$ and $\t\in[-1,1]$
$$
\pd_t \bigl( u_0^{-1}u_1 \bigr) 
\;=\; - u_0^{-1}\dot u_0 u_0^{-1}u_1 + u_0^{-1}\dot u_1 
\;=\; u_0^{-1} \ta_i(\t)^*( A_0 - A_1 ) \, u_1 .
$$
Hence
\begin{align*}
{\rm dist}_\rG\bigl( \hol_{\ta_i(\t)}(A_0) \,,\,  \hol_{\ta_i(\t)}(A_1)  \bigr)
&= {\rm dist}_\rG\bigl( \one \,,\,  u_0(1)^{-1} u_1(1) \bigr)  \\
&\leq \int_0^1 \bigl| \pd_t \bigl( u_0(t)^{-1}u_1(t) \bigr) \bigr| \,\dt \\
&\leq \int_0^1 \bigl| \ta_i(\t)^*( A_0 - A_1 ) \bigr| \,\dt .
\end{align*}
Next, for every $i=1,\ldots,g$ there exists a set $V_i\subset[-1,1]$ of measure
$|V_i|\geq 2-\frac 1g$ such that for all $\t\in V_i$
\begin{align*}
\int_0^1 \bigl| \ta_i(\t)^*( A_0 - A_1 ) \bigr| \,\dt 
&\;\leq\; g \int_{-1}^1 \int_0^1 \bigl| \ta_i(\t)^*( A_0 - A_1 ) \bigr| \,\dt\,\rd\t \\
&\;\leq\; C \int_{\ta_i} | A_0 - A_1 | \;\leq\; C \| A_0 - A_1 \|_{L^1(\S)}  .
\end{align*}
Here the constant $C$ only depends on the embeddings $\ta_i$. Now the claim (ii) is 
true for all $\t\in\bigcap_{i=1}^g V_i$, which is nonempty.
In case (i) we similarly find $\t\in[-1,1]$
such that 
$$
\int_0^1 \bigl| \ta_i(\t)^*( \pd_s A(0) ) \bigr| \,\dt \;\leq\; C  \| \pd_s A (0) \|_{L^1(\S)} .
$$
Now with $\th=\hol_{\ta_i(\t)}(A(0))$ we obtain as above
\begin{align*}
\bigl| \pd_s|_{s=0} \hol_{\ta_i(\t)}(A(s)) \bigr|_{\rT_\th\rG}
&= \lim_{s\to 0} |s|^{-1}
{\rm dist}_\rG\bigl( \hol_{\ta_i(\t)}(A(0)) \,,\,  \hol_{\ta_i(\t)}(A(s))  \bigr) \\
&\leq \lim_{s\to 0} \int_0^1 \biggl| \ta_i(\t)^*\biggl( \frac{A(0)-A(s)}s \biggr) \biggr| \,\dt \\
&= \int_0^1 \bigl| \ta_i(\t)^* \pd_s A(0) \bigr| \,\dt .
\end{align*}

\vspace{-5.5mm}

\hspace*{\fill}$\Box$
\medskip

Now consider the extension problems in lemma~\ref{lem b}. 
Given connections in $\cL_Y$, the above lemma provides a control of the holonomies 
based at some point $z(\t)$. This point can vary in a neighbourhood of $z\in\S$.
However, for any such basepoint, the sections (\ref{section}) will provide
flat connections on $Y$ with the holonomy of the given connections on $\S$.
So on $\pd Y=\S$, these connections only differ by a gauge transformation.
Thus we require the following extension construction for gauge transformations.
Here and in the following $\rd^\S_\X$ for $\X\in\cA(Y)$ denotes the exterior derivative 
on $\cA(\S)$ associated with the connection $\X|_\S$.

\begin{lem} \label{lem gauge ext}
There is a constant $C$ such that the following holds for any
connection $\X=\X_j(\Th)\in\cA_{\rm flat}(Y)$, $\Th\in B_\D(\Th^0_j)$
in the finitely many sections (\ref{section}).
\begin{enumerate}
\item 
For all $\x\in\cC^\infty(\S,\cg)$ there exists
$\tx\in\cC^\infty(Y,\cg)$ such that $\tx|_{\pd Y}=\x$ and
$$
\bigl\| \rd_\X\tx \bigr\|_{L^3(Y)}
\;\leq\; C \| \rd^\S_\X\x \|_{L^2(\S)} .
$$
\item
For all $u\in\cG(\S)$ there exists
$\tu\in\cG(Y)$ such that $\tu|_{\pd Y}=u$ and
$$
\bigl\| \tu^*\X - \X \bigr\|_{L^3(Y)}
\;\leq\; C \| u^*\X|_\S - \X|_\S \|_{L^2(\S)} .
$$
\end{enumerate}
\end{lem}

For (ii) note that a smooth map $\S\to\rG$ can always be extended to $Y\to\rG$
since by assumption $\pi_1(\rG)=0$ (so extensions to discs in $Y$ with boundary
in $\S$ exist), and for general Lie groups $\pi_2(\rG)=0$ (so these extensions
can be matched up).
We will moreover use the following quantitative result
with $N=\rG$ and thus $\ell=2$, where the Sobolev spaces of maps into $N\subset\R^k$
are understood as
$$
W^{1,q}(\Om,N) = \bigl\{ u \in W^{1,q}(\Om,\R^k) \st \forall' x\in\Om : u(x)\in N \bigr\} .
$$

\begin{thm} {\bf \cite{Hang}} \label{thm p ext}
Let $N\subset\R^k$ be a smooth connected compact Riemannian manifold 
with $\pi_i(N)=0$ for all $i=1,\ldots,\ell$.
Then the following holds for all $1<q<\ell+2$.
Let $\Om\subset\R^m$ be an open, bounded domain with piecewise smooth boundary.
Then there exists a constant $C$ such that for any $u\in W^{1-\frac 1q,q}(\pd\Om,N)$
there exists $\tu\in W^{1,q}(\Om,N)$ with $\tu|_{\pd\Om}=u$ such that
$$
\| \rd\tu \|_{L^q(\Om)} \leq C \| u \|_{W^{1-\frac 1q,q}(\pd\Om)} .
$$
In particular, if $\Om$ is simply connected and if we fix $1<\bq\leq q$ and $1<\bp\leq p$ such that
$p\geq\frac{m-1}m q$, $\bp\geq\frac{m-1}m \bq$, then there is a constant $C$ such that
for any $u\in W^{1,p}(\pd\Om,N)$ there exists $\tu\in W^{1,q}(\Om,N)$ with $\tu|_{\pd\Om}=u$ such that
$$
\| \rd\tu \|_{L^\bq(\Om)} \leq C \| \rd u \|_{L^\bp(\pd\Om)} .
$$
\end{thm}

The first part is \cite[Thm~2.1]{Hang}, and the second part is an easy
consequence:
One has $W^{1,p}(\pd\Om)\hookrightarrow W^{1-\frac 1q,q}(\pd\Om)$ 
and the trace $W^{2,\bq^*}(\Om)\hookrightarrow W^{1,\bp^*}(\pd\Om)$ 
by e.g.\ \cite[Thms~7.8,~5.22]{Adams}.
If $\Om$ is simply connected, then the operator $(\rd,\rd^*,\cdot|_{\pd\Om})$ is injective,
so as in the proof of lemma~\ref{lem gauge ext}~(i) one finds for all $\a\in\Om^1(\Om)$
$$
\| \a \|_{L^\bq(\Om)} \;\leq\;
C \bigl( \| \rd\a \|_{(W^{1,\bq^*}(\Om))^*} +  \| \rd^*\a \|_{(W^{1,\bq^*}(\Om))^*}
+  \| \a|_{\pd\Om} \|_{L^\bp(\pd\Om)} \bigr) .
$$
The proof in \cite{Hang} uses the solution $v\in W^{1,q}(\Om,\R^k)$
of $\rd^*\rd v=0$ with $v|_{\pd\Om}=u$, for which in this case
$\|\rd v\|_{L^\bq(\Om)}\leq C \|\rd u\|_{L^\bp(\pd\Om)}$. 
Variation of a 'centre' of a retraction to $N$
%
%
%
%
then gives $\tu\in W^{1,q}(\Om,N)$ with $\|\rd\tu\|_{L^q(\Om)}\leq C \|\rd v\|_{L^q(\Om)}$.
This centre $a\in\R^k$ can be found to simultaneously yield the same estimate for $\bq$.\\

\noindent
{\bf Proof of lemma \ref{lem gauge ext}: } \\
For (i) we determine $\tx\in\cC^\infty(Y,\cg)$ by solving the Dirichlet problem
\[
\rd_\X^*\rd_\X \tx = 0, \qquad\qquad
\tx|_{\pd Y} = \x .
\]
The operator $(\rd_\X^*\rd_\X,\cdot|_{\pd Y})$ on $W^{2,2}(Y,\cg)$ 
is a compact perturbation of the standard Dirichlet operator $(\laplace,\cdot|_{\pd Y})$,
so it is a Fredholm operator of index $0$.
It is surjective since its kernel equals $\ker(\rd_\X,\cdot|_{\pd Y})=\{0\}$,
where a solution of $\rd_\X\e=0$ is uniquely determined by its value at any one point 
via integration along paths.
For $\x\in\cC^\infty(\S,\cg)$ the smoothness of the
solution $\tx$ follows from elliptic regularity.

The estimate for $\rd_\X\tx$ will be provided by the following Hodge type estimate:
There exists a constant $C$, independent of $\X$, such that for all $\a\in\Om^1(Y,\cg)$.
\begin{equation} \label{1form est}
\| \a \|_{L^3(Y)} \;\leq\;
C \bigl( \| \rd_\X\a \|_{L^{\frac 32}(Y)} +  \| \rd_\X^*\a \|_{L^{\frac 32}(Y)}
+  \| \a|_{\pd Y} \|_{L^{2}(\S)} \bigr) .
\end{equation}
If we put in $\a=\rd_\X\tx$, then $\rd_\X\rd_\X\tx=0$ since $\X$ is flat and
$\rd_\X^*\rd_\X\tx=0$ by construction.
So it remains to establish (\ref{1form est}). 
If we consider the normal and tangential components of the $1$-forms on $Y$ separately,
then this estimate deals with the operator $\rd_\X^*\rd_\X$ with Dirichlet boundary conditions 
for the tangential components.
From $\rd_\X^*\a$ one also has a Neumann boundary condition for the normal component
in terms of the tangential components.
So a combination of Dirichlet estimates for the tangential components and a Neumann
estimate for the normal component will imply (\ref{1form est}).

More precisely, one can use \cite[Thm~5.3]{W} to obtain $W^{1,\frac 32}$-estimates
for $\a(X)$, where $X\in\G(\rT Y)$ is either tangential to $\pd Y$ (in which case one
uses test functions $\p\in\cC^\infty(Y,\cg)$ with $\p|_{\pd Y}=0$), or $X$ is normal to $\pd Y$ 
(and one uses test functions with $\frac\pd{\pd\n}\p|_{\pd Y} = 0$).
In both cases one then has the following estimates, where the constant $C$ depends on $\X$.
Firstly, the boundary term vanishes in
\begin{align*}
\biggl| \int_Y \la \a \,,\, \rd\cL_X\p \ra \biggr|
&= \biggl| \int_Y \la \rd_\X^*\a \,,\, \cL_X\p \ra
+ \int_Y \la *[\X \wedge *\a] \,,\, \cL_X\p \ra
+ \int_{\pd Y} \la *\a \,,\, \cL_X\p \ra \biggr| \\
&\leq C \bigl( \| \rd_\X^*\a \|_{L^{\frac 32}(Y)} 
      + \| \a \|_{(W^{1,\frac 32}(Y))^*} \bigr) \| \p \|_{W^{2,\frac 32}(Y)} .
\end{align*}
This also uses the Sobolev inequality 
$\|\p\|_{W^{1,3}(Y)}\leq C \| \p \|_{W^{2,\frac 32}(Y)}$.
Secondly, one

\noindent
 can use the Sobolev embedding $W^{2,\frac 32}(Y)\hookrightarrow W^{1,2}(\pd Y)$
(see \cite[5.22]{Adams}) to obtain
\begin{align*}
& \biggl| \int_Y \la \a \,,\, \rd^* ( i_X g\wedge\rd\p ) \ra \biggr| \\
&= \biggl| \int_Y \la \rd_\X\a \,,\, i_X g\wedge\rd\p \ra
- \int_Y \la [\X \wedge \a] \,,\, i_X g\wedge\rd\p \ra
- \int_{\pd Y} \la \a \wedge * ( i_X g\wedge\rd\p ) \ra \biggr| \\
&\leq C \bigl( \| \rd_\X\a \|_{L^{\frac 32}(Y)} + \| \a \|_{(W^{1,\frac 32}(Y))^*} 
              + \| \a|_{\pd Y} \|_{L^2(\S)} \bigr)  \| \p \|_{W^{2,\frac 32}(Y)} .
\end{align*}
These two estimates can be considered as weak Laplace equations on $\a(X)$ 
with inhomogenous Dirichlet or Neumann boundary conditions respectively.
The according estimates sum up to 
$$
\| \a \|_{L^3(Y)} \;\leq\;
C \bigl( \| \rd_\X\a \|_{L^{\frac 32}(Y)} +  \| \rd_\X^*\a \|_{L^{\frac 32}(Y)}
+  \| \a|_{\pd Y} \|_{L^{2}(\S)} +  \| \a \|_{(W^{1,\frac 32}(Y))^*} \bigr) .
$$
Finally, the last term can be dropped since the embedding 
$L^3(Y)\hookrightarrow (W^{1,\frac 32}(Y))^*$ is the dual of a compact Sobolev embedding,
and the operator $(\rd_\X,\rd_\X^*, \cdot|_{\pd Y})$ is injective.
To see the latter consider an element $\a\in\Om^1(Y,\cg)$ of the kernel.
We can write it as $\a=\rd_\X\e$ for some $\e\in\cC^\infty(Y,\cg)$ with 
$\e|_{\pd Y}=0$.
\footnote{
Since $F_\X=0$ and $\rd_\X\a=0$ this is true on simply connected subsets of $Y$.
We can moreover prescribe $\e|_{\pd Y}=0$ since $\a|_{\pd Y}=0$.
Now $Y$ can be covered with simply connected domains whose intersections are connected
and meet $\pd Y$. (The $1$-skeleton of $Y$ can be pushed to $\pd Y$.) So if $\e$ and 
$\e'$ are each determined on one of these domains, then they have to match up on the
intersection. This is since $\rd_\X(\e-\e')=0$ with $\e=\e'$ at one point only has the
trivial solution $\e=\e'$.
}
Then $\rd_\X^*\rd_\X\e=0$ with $\e|_{\pd Y}=0$ implies $\a=\rd_\X\e=0$ by partial 
integration.
Thus (\ref{1form est}) holds for every $\X\in\cA_{\rm flat}(Y)$.

The constant in (\ref{1form est}) depends continuously on $\X$ with respect 
to the $L^\infty$-norm.
It can be chosen uniform since we only consider smooth connections $\X$ 
that are parametrized by a finite number of compact sets 
$ B_\D(\Th^0_j)\subset\rG^g$.\\

In (ii) we need to extend $u:\pd Y \to \rG$ to $\tu: Y \to \rG$.
Our construction will make use of theorem~\ref{thm p ext}, where we fix an embedding
$\rG\subset\R^k$ and some $2<p<\frac 83$.
We recall the diffeomorphisms $\psi_i:[0,1]\times D \to Z_i\subset Y$ and
denote $D(\t):=\psi_i(\t,D)\subset Y$ with the orientation induced by $\psi_i$.
By construction the connection $\X$ vanishes over $D(\t)$ for all $\t\in[\frac 34,1]$.
So given any $u\in\cG(\S)$ we find $\t_i\in[\frac 34,1]$ for every $i=1,\ldots,g$ such that
$$
\int_{\pd D} |\psi_i(\t_i)^*\rd u |^2 \;\leq\; 4 \int_{[\frac 34,1]\times\pd D} |\psi_i^*\rd u |^2 
\;=\; 4 \int_{[\frac 34,1]\times\pd D} \bigl| \psi_i^*\bigl( u^*\X - \X \bigr)\bigr|^2 .
$$
Since the $\psi_i$ are fixed we then have with a uniform constant $C$ for all $i=1,\ldots,g$ 
$$
\|\rd u \|_{L^2(\pd D(\t_i))} \;\leq\; C \| u^*\X|_\S - \X|_\S \|_{L^2(\S)} .
$$
Now theorem~\ref{thm p ext} on $\Om=D(\t_i)\subset\R^2$ with $q=p>2$ as fixed and 
$\bq=\bp=2$ gives $\tu_i\in W^{1,p}(D(\t_i),\rG)$ with
$\tu_i|_{\pd D(\t_i)}= u|_{\pd D(\t_i)}$ and
$$
\|\rd\tu_i\|_{L^2(D(\t_i))} \;\leq\; C \| u^*\X|_\S - \X|_\S \|_{L^2(\S)} .
$$
Next, fix an embedding $Y\subset\R^3$ and cut $Y$ open to obtain the simply connected 
open manifold $Y_\bt={\rm int}(Y)\setminus \bigcup_{i=1}^g D(\t_i)$. 
For any choice of $\bt=(\t_i)\in[\frac 34,1]^g$ this is diffeomorphic to the standard
domain ${\rm int}(Y)\setminus\bigcup Z_i\subset\R^3$ with a uniform bound on every given derivative.
Thus we can apply theorem~\ref{thm p ext} with a uniform constant to all these domains.
Their piecewise smooth boundary then is 
$$
\pd Y_\bt = \S_\bt \cup \bigcup_{i=1}^g D(\t_i) \cup \bigcup_{i=1}^g \bar D(\t_i) 
\qquad\text{with}\qquad
\S_\bt = \S\setminus \bigcup_{i=1}^g \pd D(\t_i) .
$$
Here $D(\t_i)$ is the boundary component attached to $\psi_i([0,\t_i)\times D)\subset Y_\bt$, whereas
$\bar D(\t_i)$ with the reversed orientation is attached to $Y_\bt\setminus\psi_i([0,\t_i)\times D)$.
Now recall that $\X|_{Z_i}=v_i^{-1}\rd v_i$, where $v_i$ is smooth on $Z_i=\psi_i([0,1]\times D)$ and 
satisfies $v_i|_{\psi_i([0,\frac 12]\times D)}\equiv \one$ and $v_i|_{D(\t_i)}\equiv \th_i^{-1}\in\rG$.
So we can write $\X|_{Y_\bt}=v^{-1}\rd v$, where $v\in\cC^\infty(Y_\bt,\rG)$ is given by 
$v=v_i$ on $\psi_i([0,\t_i)\times D)$ and $v\equiv\one$ on the complement.
With this we define
$$
w := \left\{
\begin{array}{cl}
v \,u \,v^{-1} &;\text{on}\; \S_\bt \\
\tu_i  &;\text{on}\; \bar D(\t_i) \\
\th_i^{-1}\tu_i \,\th_i  &;\text{on}\; D(\t_i)
\end{array}
\right.
\qquad\in\; W^{1,p}(Y_\bt,\rG) .
$$
This gauge transformation is chosen such that on $\S_\bt$
\begin{align*}
w^{-1} \rd w 
\;=\; v\, u^{-1} v^{-1} \rd v \,u \,v^{-1} + v \,u^{-1} \rd u \,v^{-1} - v \, v^{-1} \rd v \,v^{-1} 
\;=\;  v ( u^*\X  - \X ) v^{-1} .
\end{align*}
So we can apply theorem~\ref{thm p ext} on $Y_\bt\hookrightarrow\R^3$ with $q=\frac 32 p$, 
$\bp=3$, and $\bq=2$ to obtain $\tw\in W^{1,\frac 32 p}(Y_\bt,\rG)$ such that 
$\tw|_{\pd Y_\bt}=w$ and
\begin{align*}
\| \rd\tw \|_{L^3(Y_\bt)} 
&\;\leq\; C \| \rd w \|_{L^2(\pd Y_\bt)} \\
&\;\leq\; C \bigl( \| u^*\X  - \X \|_{L^2(\S_\bt)} 
                + \|\rd\tu_i \|_{L^2(D(\t_i))} 
                + \|\th_i^{-1}\rd\tu_i\th_i \|_{L^2(D(\t_i))} \bigr) \\
&\;\leq\; C \| u^*\X|_\S - \X|_\S \|_{L^2(\S)} .
\end{align*}
Now $\tu:= v^{-1}\tw \,v \in W^{1,\frac 32}(Y_\bt,\rG)$ satisfies
$\tu|_{\S_\bt} = u|_{\S_\bt}$ and $\tu|_{D(\t_i)} = \tu_i = \tu_{\bar D(\t_i)}$,
so it matches up to $\tu\in W^{1,\frac 32}(Y,\rG)$.
Also, 
$$
\bigl( \tu^*\X - \X \bigr)\bigr|_{Y_\bt}
\;=\; v^* \tw^* (v^{-1})^* v^{-1}\rd v - v^{-1}\rd v
\;=\; v^{-1} \tw^{-1}\rd \tw \, v ,
$$
and hence
$$
\| \tu^*\X - \X \|_{L^3(Y)} 
\;=\;\| \rd\tw \|_{L^3(Y_\bt)} 
\;\leq\; C \| u^*\X|_\S - \X|_\S \|_{L^2(\S)} .
$$
Finally, we need a smooth approximation of $\tu$ that so far is only
continuous.
In case $\| u^*\X|_\S - \X|_\S \|_{L^2(\S)}=0$ we have $\rd\tu = \tu\X - \X\tu$,
where $\X$ is smooth, so automatically $\tu\in\cG(Y)$.
Otherwise we can find a smooth approximation $\bar u\in\cG(Y)$ of the map
$\tu\in W^{1,\frac 32 p}(Y,\rG)\subset\cC^0(Y,\rG)$ with fixed boundary values\footnote
{
Pick any extension $v\in\cC^\infty(Y,\rG)$ of $\tu|_{\pd Y}$.
Then $\tu-v\in W^{1,\frac 32 p}(Y,\R^k)$ has zero boundary values and thus can be
approximated by $w\in\cC^\infty_0(Y,\R^k)$. Now $v+w\in\cC^\infty(Y,\R^k)$ is $\cC^0$-close
to $\tu$ and already identic to it on $\pd Y$. So a projection from a neighbourhood of 
$\rG\subset\R^k$ to $\rG$ composed with $v+w$ yields the required approximation.
}
$\bar u|_{\pd Y}= \tu|_{\pd Y}= u$ 
and $\|\bar u - \tu\|_{W^{1,\frac 32 p}(Y,\R^k)}\leq \min ( 1 \,,\, \| u^*\X|_\S - \X|_\S \|_{L^2(\S)} )$.
This is an approximation in $W^{1,3}(Y)$ as well as $\cC^0(Y)$.
So we introduce the notation $\rd_\X\tu=\tu ( \tu^*\X - \X )=\rd\tu + \X\tu - \tu\X$ 
to estimate
\begin{align*}
&\| \bar u^*\X - \X \|_{L^3(Y)} \\
&\leq \| \tu^*\X - \X \|_{L^3(Y)} 
     + \bigl\| \bar u^{-1}\bigl(\rd_\X\bar u - \rd_\X\tu \bigr) \bigr\|_{L^3(Y)}  
     + \bigl\| \bigl(\bar u^{-1} - \tu^{-1}\bigr) \rd_\X\tu \|_{L^3(Y)} \\
&\leq  C \| u^*\X|_\S - \X|_\S \|_{L^2(\S)} .
\end{align*}
The constant $C$ again depends on $\X\in\cA^{0,3}(Y)$, but since this only varies in a compact set,
it can be chosen uniform.

\QED

\noindent
{\bf Proof of lemma \ref{lem b} (i): } \\
For a given smooth path $A : (-\ep,\ep)\to \cL_Y\cap\cA(\S)$ 
let $\Th=(\th_i)\in\cC^\infty((-\ep,\ep),\rG^g)$ be given by
$\th_i(s)=\hol_{\ta_i(\t)}(A(s))$. We pick $\t\in[-1,1]$ as in lemma~\ref{lem holonomy}~(i), so
$$
\bigl| \pd_s\Th(0) \bigr|_{\rT_{\Th(0)}\rG^g}
 \leq C  \| \pd_s A (0) \|_{L^1(\S)} .
$$
We can also pick one of the fixed $\Th^0_j\in\rG^g$ with
$\Th(s)\in B_\D(\Th^0_j)$ for all $s\in(-\ep,\ep)$ for some smaller $\ep>0$.
(Note that it suffices to construct $\tA(s) \in \cA_{\rm flat}(Y)$ for a neighbourhood of $s=0$.
Then we can arbitrarily extend it to a larger interval.)
Now we can use the chart (\ref{chart}) with $z=z(\t)$ to write
$A(s)=u(s)^*\X_j(\Th(s))|_\S$ with a smooth path $u:(-\ep,\ep)\to\cG_z(\S)$.
So we have
$$
\pd_s A (0) = 
u(0)^{-1} \bigl( \rT_{\X_0} \X_j (\pd_s\Th(0)) \bigr|_\S + \rd^\S_{\X_0} \x \bigr) \, u(0)
$$
with $\x = \pd_s u(0) u(0)^{-1} \in \cC^\infty_z(\S,\cg)$ and
$\X_0=\X_j(\Th(0))$. 
Here the operators $\rT_{\X_0}\X_j |_\S : \rT B_\D(\Th^0_j) \to L^2(\S,\rT^*\S\otimes\cg)$ 
bounded, and we can choose a uniform constant on all of $B_\D(\Th^0_j)$ 
for all $j=1,\ldots,N$. So we have with another uniform constant $C$
$$
\bigl\| \rd^\S_{\X_0} \x \bigr\|_{L^2{\S}}
\;\leq\; \bigl\| \pd_s A (0) \bigr\|_{L^2{\S}} 
       + \bigl\| \rT_{\X_0} \X_j |_\S \bigr\| \bigl|\pd_s\Th(0)\bigr|_{\rT_{\Th(0)}\rG^g}
\;\leq\; C \bigl\| \pd_s A (0) \bigr\|_{L^2{\S}}  .
$$
Next, lemma~\ref{lem gauge ext} provides $\tu_0\in\cG(Y)$ with $\tu_0|_{\pd Y}=u$
and $\tx\in\cC^\infty(Y,\cg)$ such that $\tx|_{\pd Y}=\x$ and
$$
\|\rd_{\X_0}\tx\|_{L^3(Y)}
\;\leq\; C \|\rd^\S_{\X_0}\x\|_{L^2(\S)} .
$$
If we now define $\tA : (-\ep,\ep)\to \cA_{\rm flat}(Y)$ by 
$\tA(s)=( \exp(s\tx)\, \tu_0 )^* \X_j(\Th(s))$,
then indeed $\pd_s\tA(0)|_{\pd Y}=\pd_s A(0)$ and
\begin{align*}
\|\pd_s\tA(0)\|_{L^3(Y)} 
&\;=\; \bigl\| \tu_0^{-1}\bigl(\rT_{\X_0} \X_j (\pd_s\Th(0)) 
                                           + \rd_{\X_0} \tx \bigr)\, \tu_0 \bigr\|_{L^3(Y)}  \\
&\;\leq\; \bigl\| \rT_{\X_0} \X_j \bigr\|  | \pd_s\Th(0) | + \|\rd_{\X_0} \tx \|_{L^3(Y)} 
\;\leq\; C \|\pd_s A(0)\|_{L^2(\S)} .
\end{align*}
Here we have moreover chosen a uniform constant of continuity for the operators
$\rT_{\X_0}\X_j : \rT B_\D(\Th^0_j) \to L^3(Y,\rT^*Y\otimes\cg)$ 
on the compact domains $B_\D(\Th^0_j)$ for all $j=1,\ldots,N$. 
\QED

\noindent
{\bf Proof of lemma \ref{lem b} (ii): } \\
Let $A_0, A_1\in\cL_Y\cap\cA(\S)$ be given.
We will prove the lemma by construction, assuming that $A_0 = \X_j(\P^0)|_\S$ for some 
$\P^0=(\p^0_i)\in B_{\frac 12 \D}(\Th^0_j)$.

In general, we have $u_0\in\cG(\S)$ such that $A_0 = u_0^*\X_j(\P^0)|_\S$.
The construction below then gives extensions 
$\tA_0, \tA_1\in\cA_{\rm flat}(Y)$ of $(u_0^{-1})^*A_0$ and $(u_0^{-1})^*A_1$.
Moreover, lemma~\ref{lem holonomy} provides $\tu_0\in\cG(Y)$ such that 
$\tu_0|_{\pd Y}=u_0$. Then $\tu_0^*\tA_0$ and $\tu_0^*\tA_1$ are extensions
of $A_0$ and $A_1$, and the estimate on $\tA_0-\tA_1$ also yields
\begin{align*}
\bigl\|\tu_0^*\tA_0-\tu_0^*\tA_1\bigr\|_{L^3(Y)}
&= \|\tA_0-\tA_1\|_{L^3(Y)}  \\
&\leq C_Y \bigl\|(u_0^{-1})^*A_0-(u_0^{-1})^*A_1\bigr\|_{L^2(\S)}
= C_Y \|A_0-A_1\|_{L^2(\S)} .
\end{align*}
So from now on suppose that $A_0 = \X_j(\P^0)|_\S$.
Then we already have the extension $\tA_0:=\X_j(\P^0)\in\cA_{\rm flat}(Y)$.
Note that $\hol_{\ta_i(\t)}(A_0)=\p^0_i$ for all $\t\in[-1,1]$.
Lemma~\ref{lem holonomy}~(ii) then provides $\t\in[-1,1]$ such that for all $i=1,\ldots,g$
$$
{\rm dist}_{\rG}\bigl( \p^0_i \,,\,  \hol_{\ta_i(\t)}(A_1)  \bigr) 
\;\leq\; C \| A_0 - A_1 \|_{L^1(\S)}. 
$$
If $\| A_0 - A_1 \|_{L^1(\S)}\leq \frac\D{2C}$ then this implies
$\P:=(\hol_{\ta_i(\t)}(A_1))_{i=1,\ldots,g}\in B_{\D}(\Th^0_j)$.
In that case we have found a flat connection $\tA:=\X_j(\P)$ on $Y$ whose
holonomies (based at $z(\t)$) coincide with those of $A_1$, and
\begin{align}
&\|\tA_0-\tA\|_{L^3(Y)} + \| A_0-\tA|_\S\|_{L^2(\S)}  \nonumber\\
&\;=\; \bigl\| \X_j(\P_0)-\X_j(\P) \bigr\|_{L^3(Y)} 
     + \bigl\| \bigl( \X_j(\P_0)-\X_j(\P) \bigr)|_\S \bigr\|_{L^2(Y)} \nonumber\\
&\;\leq\; C \,{\rm dist}_{\rG^g}\bigl( \P_0\,,\,\P \bigr)
\;\leq\; C \|A_0-A_1\|_{L^1(\S)} .  \label{tA tA0}
\end{align}
Here and in the following, all uniform constants are denoted by $C$.
We have in particular used the fact that the sections $\X_j$ and $\X_j|_\S$ are 
smooth on a compact set, so they are Lipschitz continuous with uniform constants.

In case $\| A_0 - A_1 \|_{L^1(\S)}\geq \frac\D{2C}$ we also use the
sections (\ref{section}) to find a flat connection $\tA:=\X_{j'}(\P)$ on $Y$ 
with the same holonomies (based at $z(\t)$) as $A_1$. The sections
are uniformly bounded in $L^3(Y)$ and $L^2(\S)$ since they are smooth over 
a union of compact sets.
Hence there is a uniform constant $\bar C$ such that
$\|\tA_0-\tA\|_{L^3(Y)} + \| A_0-\tA|_\S\|_{L^2(\S)} \leq \bar C$, 
and thus (\ref{tA tA0}) again holds with $C=\tfrac{2C\bar C}{\D}$.

For the two flat connections $A_1$ and $\tA|_\S$ with coinciding holonomies
one then finds a gauge transformation $u\in\cG(\S)$ such that $u^*\tA|_\S = A_1$.
Now by lemma~\ref{lem gauge ext}~(ii) there exists an extension
$\tu\in\cG(Y)$ with $\tu|_\S = u$ and such that
\begin{align*}
\bigl\| \tu^*\tA - \tA \bigr\|_{L^3(Y)}
&\;\leq\; C \| u^*\tA|_\S - \tA|_\S \|_{L^2(\S)} \\
&\;\leq\; C \bigl( \| A_1 - A_0 \|_{L^2(\S)} + \| A_0 - \tA|_\S \|_{L^2(\S)} \bigr) \\
&\;\leq\; C \| A_1 - A_0 \|_{L^2(\S)} .
\end{align*}
So if we put $\tA_1:=\tu^*\tA\in\cA_{\rm flat}(Y)$, then indeed $\tA_1|_{\pd Y}=A_1$
and
\begin{align*}
\bigl\| \tA_1 - \tA_0 \bigr\|_{L^3(Y)}
\;\leq\; \bigl\| \tu^*\tA - \tA \bigr\|_{L^3(Y)}
        + \bigl\| \tA - \tA_0 \bigr\|_{L^3(Y)}  
\;\leq\; C \| A_1 - A_0 \|_{L^2(\S)} .
\end{align*}
\QED

\section{Isoperimetric inequalities}
\label{isoperimetric}

The aim of this section is to firstly introduce the local Chern-Simons functional 
and prove the isoperimetric inequality, lemma~\ref{lem c}.
Secondly, we will show how the Chern-Simons functional is related to the energy
of solutions of the boundary value problem (\ref{bvp}).
This relation will yield a control of the energy that will be the key to the
removal of singularities in the next section.

The usual Chern-Simons functional on a closed $3$-manifold $M$ is 
\begin{align*}
\CS(\X) = \half \int_M \la \X \wedge \bigl( F_\X - \tfrac 16 [\X\wedge\X] \bigr) \ra  
\qquad\qquad\forall\,\X\in\cA(M) .
\end{align*}
It is not gauge invariant, but its change 
$\CS(\X) - \CS(u^*\X) = 4\pi^2 {\rm deg}(u)$ is determined by the
degree of the gauge transformation $u\in\cG(M)$.
The negative gradient flow lines of $\CS$ are the anti-self-dual connections on $\R\times M$.
This can be seen from the fact that the differential $\rd_\X\CS:\Om^1(M;\cg)\to\R$
is given by $\a \;\mapsto\; \tint_M \la \a \wedge F_\X \ra$.

If $M$ is a manifold with boundary, then this $1$-form is not closed -- its differential is
the natural symplectic structure on $\Om^1(\pd M;\cg)$, c.f.\ \cite{Sa1}.
So it is natural to impose Lagrangian boundary conditions $\X|_{\pd M}\in\cL$. On this 
subset of connections, the above $1$-form is closed. However it is only the differential
of a multi-valued functional.
If the Lagrangian is $\cL_Y$, given by the flat connections on a handle body $Y$
restricted to the boundary $\pd Y=\S$, then this multi-valued Chern-Simons functional
can be represented as follows.
Given $\X\in\cA(M)$ with $\X|_{\pd M}\in\cL_Y$ one can find $\tX\in\cA_{\rm flat}(Y)$
with $\tX|_{\pd Y}=\X|_{\pd M}$ and use this to define
$$
\CS_{\cL_Y}(\X) = 
\half \int_M \la \X \wedge \bigl( F_\X - \tfrac 16 [\X\wedge\X] \bigr) \ra  
+\tfrac 1{12} \int_Y \la \tX \wedge [\tX\wedge\tX] \ra  .
$$
This is the actual Chern-Simons functional on the closed manifold $M\cup_\S \bar Y$
(where $\bar Y$ has the reversed orientation) of the connection given by $\X$ and
$\tX$ on the two parts. It is welldefined only up to multiples of $4\pi^2$ due to the
choice of different extensions $\tX$ of $\X|_{\pd M}$. A change of this extension
corresponds to the action of a gauge transformation on $M\cup_\S \bar Y$ that
is trivial on $M$.
(The gauge equivalence class of a flat connection on $Y$ is fixed by its
holonomies on $\pd Y$.)

Our energy identities below deal with connections $\X\in\cA([0,\pi]\times\S)$ 
with boundary values $\X|_{\p=0},\X|_{\p=\pi}\in\cL_Y$.
These can be put into the special gauge $\X=A$ with $A:[0,\pi]\to\cA(\S)$. 
So equivalently to $\CS_{\cL_Y}(\X)$, we can define the
{\bf local Chern-Simons functional} for smooth paths $A:[0,\pi]\to\cA(\S)$ 
with endpoints $A(0),A(\pi)\in\cL_Y$ 
(that will actually be welldefined for short paths):
\begin{align}
\CS(A)
&= -\half \int_0^\pi \int_\S \la A \wedge \pd_\p A \ra \,\dph   \label{CS}\\
&\quad - \tfrac 1{12} \int_Y \la \tA(0) \wedge  [\tA(0)\wedge\tA(0)] \ra  
       + \tfrac 1{12} \int_Y \la \tA(\pi) \wedge [\tA(\pi)\wedge\tA(\pi)] \ra   , \nonumber
\end{align}
where $\tA(0),\tA(\pi)\in\cA_{\rm flat}(Y)$ such that $\tA(0)|_{\pd Y}=A(0)$,
$\tA(\pi)|_{\pd Y}=A(\pi)$, and
\begin{equation} \label{ext cond}
\|\tA(0)-\tA(\pi)\|_{L^3(Y)} \leq C_Y \|A(0)-A(\pi)\|_{L^2(\S)} .
\end{equation}
Here $C_Y$ is the constant from lemma~\ref{lem b}~(ii), which ensures the existence
of the extensions $\tA(0)$ and $\tA(\pi)$.
This $\CS(A)$ equals the above $\CS_{\cL_Y}(\X)$ in the special gauge.
So a priori it is defined only up to multiples of $4\pi^2$ due to the
freedom in the choice of the extensions $\tA(0),\tA(\pi)$.
However, we will see below that for sufficiently 
small $\int_0^\pi \|\pd_\p A\|_{L^2(\S)}$ this Chern-Simons functional is welldefined,
i.e.\ any choice of extensions $\tA(\pi), \tA(0)$ that satisfies (\ref{ext cond})
will give the same value for $\CS(A)$.\\

\noindent
{\bf Proof of lemma~\ref{lem c} :}\;
Let $A:[0,\pi]\to\cA(\S)$ be a smooth path with $A(0),A(\pi)\in\cL_Y$
and $\int_0^\pi \|\pd_\p A\|_{L^2(\S)}\leq\ep$, where $\ep>0$ will be fixed
later on.
Consider any flat connections $\tA(0),\tA(\pi)\in\cA_{\rm flat}(Y)$ 
such that $\tA(0)|_{\pd Y}=A(0)$, $\tA(\pi)|_{\pd Y}=A(\pi)$, and
(\ref{ext cond}) holds.
With these we calculate
\begin{align*}
&\int_0^\pi \int_\S \la A \wedge \pd_\p A \ra \, \dph \\
&= \int_0^\pi\int_\S \la \bigl( A(0) + \tint_0^\p \pd_\p A(\th) \,\rd\th \bigr)
                        \wedge \pd_\p A(\p) \ra \, \dph  \\
&= \int_0^\pi\int_0^\p \int_\S \la \pd_\p A(\th) \wedge \pd_\p A(\p) \ra \,\rd\th\,\dph 
  + \int_\S \la A(0) \wedge  A(\pi) \ra  - \la A(0) \wedge  A(0) \ra  \\
&= \int_0^\pi \hspace{0pt} \int_0^\p \hspace{0pt} \int_\S 
                                 \la \pd_\p A(\th) \wedge \pd_\p A(\p) \ra \,\rd\th\,\dph 
 + \int_Y \la \rd \tA(0) \wedge  \tA(\pi) \ra - \la \tA(0) \wedge  \rd \tA(\pi) \ra
\end{align*}
Now use the fact that $F_{\tA(0)}=F_{\tA(\pi)}=0$ 
and choose $\ep\leq \frac{6}{C_Y^3}$ to obtain
\begin{align*}
\CS(A)
&= -\half \int_0^\pi\int_0^\p \int_\S \la \pd_\p A(\th) \wedge \pd_\p A(\p) \ra \,\rd\th \,\dph \\
&\quad + \tfrac 14 \int_Y \la [\tA(0)\wedge \tA(0)] \wedge  \tA(\pi) \ra
                        - \la \tA(0) \wedge [\tA(\pi)\wedge\tA(\pi)] \ra \\
&\quad - \tfrac 1{12} \int_Y \la [\tA(0)\wedge\tA(0)] \wedge \tA(0) \ra  
                           - \la [\tA(\pi)\wedge\tA(\pi)] \wedge \tA(\pi) \ra  \\
&= -\half \int_0^\pi\int_0^\p \int_\S \la \pd_\p A(\th) \wedge \pd_\p A(\p) \ra \,\rd\th \,\dph \\
&\quad - \tfrac 1{12} \int_Y \la \bigl[ (\tA(0)-\tA(\pi)) \wedge (\tA(0)-\tA(\pi)) \bigr] 
                  \wedge \bigl(\tA(0)-\tA(\pi)\bigr) \ra    \\
\Rightarrow \bigl| \CS(A) \bigr|
&\leq \half \left( \int_0^\pi \bigl\| \pd_\p A \bigr\|_{L^2(\S)} \,\dph \right)^2
 + \tfrac 1{12}  \left( \bigl\| \tA(0)-\tA(\pi) \bigr\|_{L^3(Y)} \right)^3  \\
&\leq \bigl( \half + \tfrac {C_Y^3}{12} \bigl\| A(0)-A(\pi) \bigr\|_{L^2(\S)} \bigr)
 \left( \int_0^\pi \bigl\| \pd_\p A \bigr\|_{L^2(\S)} \,\dph \right)^2 \\
&\leq  \left( \int_0^\pi \bigl\| \pd_\p A \bigr\|_{L^2(\S)} \,\dph \right)^2 
\;\leq\;\ep^2.
\end{align*}
If we choose $\ep>0$ small enough, then this implies that our choice of extensions 
will always yield values $\CS(A)\in [-\pi^2 , \pi^2]$.
As seen before, $\CS(A)$ is the usual Chern-Simons functional on the closed $3$-manifold
$\bar Y \cup_{\{\pi\}\times\S} [0,\pi] \times \S \cup_{\{0\}\times\S} Y $
of the connection given by $\tA(\pi)$, $A$, and $\tA(0)$ on the different
parts.
If we change the extensions $\tA(0)$ and $\tA(\pi)$, then this corresponds
to changing the connection on the closed manifold by one gauge transformation 
(that is nontrivial only in the interior of $Y$ and $\bar Y$). 
Hence the Chern-Simons functional will change by a multiple 
(the degree of the gauge transformation) of $4\pi^2$.
This cannot lead to another value in the interval $[-\pi^2 , \pi^2]$, hence the
value of $\CS(A)$ is uniquely determined by the condition (\ref{ext cond}) on the
extensions.
\QED

The Chern-Simons functional is the starting point for the removal of singularities
in theorem~\ref{thm B} and remark~\ref{rmk C}. In both cases, the energy on a
neighbourhood of the singularity can be expressed by the Chern-Simons functional
(of the connection on the boundary of this neighbourhood in a certain gauge).
This will yield a control on the energy near the singularity.
In the interior case, remark~\ref{rmk C}, we fix the radius $r_0>0$ 
and a metric of normal type on $B\times\S$.
We use the following notation for circles and punctured balls centered at $0$, 
$$
S_r:=\pd B_r, \qquad\quad
B^*_r:= B_r(0)\setminus\{0\} \,\subset\,\R^2, \qquad\quad
B^* := B^*_{r_0}.
$$
We then consider a connection $\X\in\cA(B^*\times\S)$ that is anti-self-dual,
\begin{equation}\label{bvp i}
*F_\X + F_\X = 0 .
\end{equation}
Using polar coordinates $r\in(0,r_0]$, $\p\in[0,2\pi]$ on $B^*$ we assume
as in remark~\ref{rmk C} that the connection is in the gauge $\X=A+R\dr$
with no $\dph$-component and $A:D\to\Om^1(\Si,\cg)$, $R:D\to\Om^0(\Si,\cg)$.
Then (\ref{bvp i}) then identifies the curvature components
$$
*F_A \;=\; r^{-1} \pd_\p R \,, \qquad\quad
*\bigr( \pd_r A - \rd_A R \bigr) \;=\; r^{-1} \pd_\p A .
$$
Hence for $0<\r\leq r_0$ the energy of the connection on $B_\r^*\times\S$ is
\begin{align}\label{energy i}
\E(\r) & := \half \int_{B^*_\r\times\S} |F_\X|^2  
= \int_0^\r \int_0^{2\pi} \bigl( \|F_A\|_{L^2(\S)}^2 + r^{-2}\|\pd_\p A\|_{L^2(\S)}^2 \bigr) \,r\,\dph\,\dr .
\end{align}
We shall see in lemma~\ref{lem isoper}~(i) that in this gauge 
the Chern-Simons functional on $S_r\times\S$ equals the energy $\E(r)$, 
which leads to a decay estimate for the energy. 

In the boundary case, theorem~\ref{thm B}, we fix a radius $r_0>0$ and a 
metric of normal type on $D\times\S$, and we denote the punctured half balls by
$$
D^*_r:= B_r(0)\setminus\{0\} \,\cap\, \H^2 , \qquad\qquad
D^* := D^*_{r_0}.
$$
We consider a connection $\X\in\cA(D^*\times\S)$ that solves the boundary value problem
\begin{equation}\label{bvp b}
\left\{\begin{array}{l}
*F_\X + F_\X = 0,\\
\X|_{(s,0)\times\Si} \in\cL_Y \quad\forall s\in[-r_0,0)\cup (0,r_0] .
\end{array}\right.
\end{equation}
Using polar coordinates $r\in(0,r_0]$, $\p\in[0,\pi]$ on $D^*$ we can always
choose a gauge $\X=A+R\dr$ with no $\dph$-component.
Then the energy function is
\begin{align}\label{energy b}
\E(\r) &:= \half \int_{D^*_\r\times\S} |F_\X|^2  
&= \int_0^\r \int_0^\pi \bigl( \|F_A\|_{L^2(\S)}^2 + r^{-2}\|\pd_\p A\|_{L^2(\S)}^2 \bigr) \,r\,\dph\,\dr .
\end{align}
We shall see that for sufficiently small $\r>0$ this energy equals the local Chern-Simons
functional $\CS(A(\r,\cdot))$, and this yields a decay estimate for the energy.

\begin{lem} \hspace{1mm} \label{lem isoper} \hspace{1mm}\\
\vspace{-5mm} 
\begin{enumerate}
\item 
Let $\X\in\cA(B^*\times\S)$ satisfy (\ref{bvp i}) and $\E(r_0)<\infty$, and 
suppose that it is in the gauge $\X=A+R\dr$ with $\P\equiv 0$.
Then for all $r\leq r_0$
$$
\E(r) \;=\; - \CS(\X|_{S_r\times\S}) 
\;\leq\; \half \left( \int_0^{2\pi} \bigl\| \pd_\p A(r,\p)  \bigr\|_{L^2(\S)} \, \dph  \right)^2 
\;\leq\; \pi r \,\dot\E(r)
$$
and hence $\E(r)\leq C r^{2\b}$ with $\b=\frac 1{2\pi}>0$ and some constant $C$.
\item
Let $\X\in\cA(D^*\times\S)$ satisfy (\ref{bvp b}) and $\E(r_0)<\infty$,
and suppose that it is in the gauge $\X=A+R\dr$ with $\P\equiv 0$.
Then there exists $0<r_1\leq r_0$ such that for all $r\leq r_1$
$$
\E(r) \;=\; - \CS(A(r,\cdot)) 
\;\leq\; \left( \int_0^\pi \bigl\| \pd_\p A(r,\p)  \bigr\|_{L^2(\S)} \, \dph  \right)^2 
\;\leq\; \pi r \,\dot\E(r)
$$
and hence $\E(r)\leq C r^{2\b}$ with $\b=\frac 1{2\pi}>0$ and some constant $C$.
\end{enumerate}
\end{lem}

Note that for every connection on $B^*\times\S$ (and similarly for $D^*\times\S$) 
with finite energy the decay of the energy $\E(r)\to 0$ as $r\to 0$ is automatic: 
The assumption $\E(r_0)<\infty$ just means that the limit
$\half\int_{(B_{r_0}\setminus B_r)\times\S} |F_\X|^2=\E(r_0)-\E(r) \to \E(r_0)$ exists 
as $r\to 0$.
Now this lemma allows to control the rate of decay of $\E(r)$ for anti-self-dual
connections or solutions of the boundary value problem (\ref{bvp b}).

The proof of lemma~\ref{lem isoper} will make use of lemma~\ref{lem mean curv}, 
which implies that 
$$
\int r^2 \|F_\X(r,\p)\|_{L^2(\S)}^2 \,\dph \;\leq\; C \, \E(2r) 
\;\underset{r\to 0}\longrightarrow\; 0 .
$$
For any smooth connection with finite energy there always exists a sequence
$r_i\to 0$ for which the above integral converges to zero.
This suffices for the proof of lemma~\ref{lem isoper}~(i), but in case (ii)
we need this control for all sufficiently small $r>0$ in order to be able to
use the local Chern-Simons functional.
Lemma~\ref{lem isoper} will only be used for the proof of theorem~\ref{thm B} and
remark~\ref{rmk C} and does not affect the other results in section~\ref{singularity},
so we can indeed use lemma~\ref{lem mean curv} in its proof.\\

\noindent
{\bf Proof of lemma~\ref{lem isoper}:}\;
We start with the interior case (i).
Let $0<\r\leq r_0$, then by assumption $\E(\r)\leq\E(r_0)$ is finite, 
i.e.\ it exists as the limit
$$
\E(\r) \;=\; \lim_{\d\to 0} \half \int_{( B_\r\setminus B_\d ) \times \S} |F_\X|^2 .
$$
Due to the anti-self-duality of $F_\X$ we can rewrite
\begin{align*}
\half \int_{( B_\r\setminus B_\d ) \times \S} |F_\X|^2
&= - \half \int_{( B_\r\setminus B_\d ) \times \S} \la F_\X \wedge F_\X \ra \\
&= - \half \int_{( B_\r\setminus B_\d ) \times \S}
     \rd \,\la \X \wedge \bigl( F_\X - \tfrac 16 [\X\wedge\X] \bigr) \ra  \\
&= - \CS(\X|_{S_\r\times\S}) + \CS(\X|_{S_\d\times\S}) .
\end{align*}
Here the Chern-Simons functional on $S_r\times\S$ for $r=\r$ and $r=\d$
is not gauge invariant but changes by multiples of $4\pi^2$ under gauge
transformations of nonzero degree. However, the special gauge 
$\X|_{S_r\times\S}=A(r,\cdot) : [0,2\pi]\to\cA(\S)$ fixes these values, 
and we obtain
\begin{align*}
\CS(\X|_{S_r\times\S})
& = -\half \int_0^{2\pi}\int_\S \la  A \wedge \pd_\p A \ra \, \dph \\
&= -\half \int_0^{2\pi}\int_\S \la \Bigl( A(r,0) + \int_0^\p \pd_\p A(r,\th) \,\rd\th \Bigr) 
                             \wedge \pd_\p A(r,\p) \ra \, \dph  \\
&= -\half \int_0^{2\pi}\int_0^\p \int_\S \la \pd_\p A(r,\th) \wedge \pd_\p A(r,\p) \ra \, \rd\th \, \dph .
\end{align*}
Hence for all $0<r\leq r_0$
\begin{align*}
2 \, \bigl| \CS(\X|_{S_r\times\S}) \bigr|
&\leq \left( \int_0^{2\pi} \bigl\| \pd_\p A(r,\p)  \bigr\|_{L^2(\S)} \, \dph  \right)^2 \\
&\leq \pi \int_0^{2\pi} \bigl\| \pd_\p A(\r,\p)  \bigr\|_{L^2(\S)}^2 \, \dph  
\;\leq\; \half \pi \int_0^{2\pi} \r^2 \bigl\| F_\X(\r,\p) \bigr\|_{L^2(\S)}^2 \,\dph  .
\end{align*}
Now we know by lemma~\ref{lem mean curv} that the last expression (and thus also the length of the path 
$A(r,\cdot)\in\cA^{0,2}(\S)$\,) goes to zero as $r\to 0$.
Thus we obtain
$$
\E(\r) 
\;=\; - \CS(\X|_{S_\r\times\S})
\;\leq\; \half \pi \int_0^{2\pi} \r^2 \bigl\| F_\X(\r,\p) \bigr\|_{L^2(\S)}^2 \,\dph  
\;\;=\; \pi \r \, \dot \E(\r) ,
$$
\begin{equation} \label{E decay}
\Rightarrow \;
\ln\E(r) \;\leq\; \ln\E(r_0) - \int_r^{r_0} (\pi \r)^{-1} \rd \r
\;=\; \ln\E(r_0) - \tfrac 1\pi \ln r_0  + \tfrac 1\pi \ln r .
\end{equation}
Hence we have $\E(r) \leq C r^{2\b} $ with $\b=\frac 1{2\pi}>0$, which proves (i).

In (ii) we also have for all $0<\r\leq r_1$ (where $r_1>0$ will be 
fixed later on)
$$
\E(\r) \;=\; \lim_{\d\to 0} \half \int_{( D_\r\setminus D_\d ) \times \S} |F_\X|^2 .
$$
We aim to express this as the difference of a functional at $r=\r$ and at $r=\d$. 
The straightforward approach as in (i) would pick up additional boundary terms
on $\{\p=0\}$ and $\{\p=\pi\}$. We eliminate these by glueing $Y$ to $\S=\pd Y$ and 
extending the connections $A(r,0),A(r,\pi)\in\cL_Y$ to flat connections on $Y$.
More precisely, the oriented boundary of $(D_\r\setminus D_\d)\times\S$ consists of 
$\{r=\r\}\cong [0,\pi]\times\S$ and $\{r=\d\}\cong [0,\pi]\times\bar\S$ and the
additional parts $\{\p=0\}\cong[\d,\r]\times\S$
and $\{\p=\pi\}\cong[\d,\r]\times\bar\S$ (where $\bar\S$ has the reversed orientation).
So we glue in $[\d,\r]\times Y$ and $[\d,\r]\times\bar Y$ to obtain the smooth $4$-manifold
$$
X(\d,\r)= [\d,\r]\times \bar Y \cup_{\{\p=\pi\}} ( D_\r\setminus D_\d )\times \S 
           \cup_{\{\p=0\}} [\d,\r]\times Y 
$$
which has the boundary component
$\bar Y \cup_{\{\p=\pi\}} [0,\pi] \times \S \cup_{\{\p=0\}} \times Y $ at $r=\r$
and with reversed orientation at $r=\d$.

Next, $A(\cdot,0)$ and $A(\cdot,\pi)$ are smooth paths in $\cL_Y\cap\cA(\S)$. 
So we can pick smooth paths of extensions 
$\tA(\cdot,0), \tA(\cdot,\pi) : [\d,\r] \to \cA_{\rm flat}(Y)$.
We also extend the functions $R|_{\p=0}$ and $R|_{\p=\pi}$ 
from $[\d,\r]\times\S$ to smooth functions $\tR_0$ and $\tR_\pi$ on $[\d,\r]\times Y$.
These extensions match up to a $W^{1,\infty}$-connection on $X(\d,\r)$,
$$
\tX = \left\{\begin{array}{ll}
\tA(\cdot,\pi) + \tR_\pi\dr &;\text{on}\; [\d,\r]\times \bar Y , \\
A + R\dr &;\text{on}\; ( D_\r\setminus D_\d ) \times\S , \\
\tA(\cdot,0) + \tR_0\dr &;\text{on}\; [\d,\r]\times Y .
\end{array}
 \right.
$$
We will choose the two paths of extensions $\tA(\cdot,0)$ and $\tA(\cdot,\pi)$ 
such that for all $\d\leq r\leq\r$ the functional $\cC(A(r,\cdot),\tA(r,0),\tA(r,\pi))$
given by (\ref{CS}) with these extensions equals the local Chern-Simons functional
$\CS(A(r,\cdot))$.
For this purpose let $\bep>0$ be the constant from lemma~\ref{lem mean curv} and choose 
$0<r_1\leq \half r_0$ such that $\E(2r_1)\leq\bep$. Then for all $0<r\leq r_1$
\begin{align*}
\biggl(\int_0^\pi \bigl\| \pd_\p A(r,\p)  \bigr\|_{L^2(\S)} \, \dph \biggr)^2
&\leq \pi \int_0^\pi \bigl\| \pd_\p A(r,\p)  \bigr\|_{L^2(\S)}^2 \, \dph \\
&\leq \tfrac \pi 2 \int_0^\pi r^2 \bigl\| F_\X(r,\p) \bigr\|_{L^2(\S)}^2 \,\dph  
\;\leq\; C \E(2r) .
\end{align*}
Now choose $r_1>0$ even smaller such that $C\E(2r_1)\leq\min(\pi^2,\ep^2)$
with $\ep>0$ from lemma~\ref{lem c}.
Then the lemma applies to $A(r,\cdot)$ for all $0<r\leq r_1$.
In particular, since $\r\leq r_1$, we can choose the two paths of extensions 
to end at $\tA(\r,0)$ and $\tA(\r,\pi)$, and hence
$\cC(A(\r,\cdot),\tA(\r,0),\tA(\r,\pi))=\CS(A(\r,\cdot))$.

Moreover we know that for all $r\in [\d,\r]$ the path $A(r,\cdot)$ is sufficiently
small for the local Chern-Simons functional $\CS(A(r,\cdot))$ to be defined
and take values in $[-\pi^2,\pi^2]$.
Now $\cC(A(r,\cdot),\tA(r,0),\tA(r,\pi))$ is a smooth function of $r\in[\d,\r]$
whose values might differ from $\CS(A(r,\cdot))$ by multiples of $4\pi^2$.
We have equality at $r=\r$ and hence by continuity for all $r\in[\d,\r]$ as claimed.
Thus we actually obtain the local Chern-Simons functional from $\CS(\tX)$ on $\pd X(\d,\r)$,
\begin{align} 
\half \int_{( D_\r\setminus D_\d ) \times \S} |F_\X|^2 
&= - \half \int_{X(\d,\r)} \la F_\tX \wedge F_\tX \ra \nonumber\\
&= - \half \int_{\pd X(\d,\r)} 
     \la \tX \wedge \bigl( F_\tX - \tfrac 16 [\tX\wedge\tX] \bigr) \ra  \nonumber\\
&= - \CS(A(\r,\cdot)) + \CS(A(\d,\cdot)) ,  \label{energy CS}
\end{align}
Here we have $F_\tX\wedge F_\tX = -|F_\X|^2 \dvol$ on $(D_\r\setminus D_\d ) \times \S$
and $F_\tX\wedge F_\tX = 0$ on $[\d,\r]\times Y$
since $F_\tX$ vanishes on the $3$-dimensional slices $\{r\}\times Y$.
Now by lemma~\ref{lem c} 
$$
\bigl| \CS(A(r,\cdot)) \bigr|
\;\leq\;  \left( \int_0^\pi \bigl\| \pd_\p A(r,\p) \bigr\|_{L^2(\S)} \,\dph \right)^2 
\;\leq\; \frac \pi 2 \int_0^\pi r^2 \bigl\| F_\X(r,\p) \bigr\|_{L^2(\S)}^2 \,\dph  .
$$
As $r\to 0$ this expression converges to zero by lemma~\ref{lem mean curv}.
Thus for all $0<\r\leq r_1$ 
$$
\E(\r) 
\;=\; - \CS(A(\r,\cdot))
\;\leq\; \frac \pi 2 \int_0^\pi \r^2 \bigl\| F_\X(\r,\p) \bigr\|_{L^2(\S)}^2 \,\dph  
\;\;=\; \pi \r \, \dot \E(\r) .
$$
As in (\ref{E decay}) this implies $\E(r) \leq C r^{2\b} $ for all $0<r\leq r_1$ 
with $\b=\frac 1{2\pi}>0$.
\QED

\section{Removal of singularities}
\label{singularity}

This section gives the proofs of theorem~\ref{thm B} and remark~\ref{rmk C}.
We will also prove a more general removable singularity result, theorem~\ref{thm rem}, 
that does not require the connections to solve an equation but only assumes 
a decay condition on the curvature. 
For solutions of (\ref{bvp}), as a consequence of the isoperimetric and by the lemma below, 
this decay condition is equivalent to the connection having finite energy.
In the case of interior singularities of anti-self-dual connections
the same is true if we assume the existence of a special gauge as in remark~\ref{rmk C}.
Throughout this section we fix metrics of normal type on $D\times\S$ and $B\times\S$.

\begin{lem}  \label{lem decay}
Let $\X$ be a smooth connection on $D^*\times\S$ or $B^*\times\S$.
Suppose that it satisfies (\ref{bvp b}) or (\ref{bvp i}) respectively.
Then the following are equivalent:
\begin{enumerate}
\item
$\E(r)\leq C r^{2\b}$ for all $r\leq r_0$ and some constants $C$ and $\b>0$.
\item
$\sup_\p \|F_\X(r,\p)\|_{L^2(\S)}\leq C r^{\b-1}$ for all $r\leq r_0$  
and constants $C$ and $\b>0$.
\item
$\|F_\X\|_{L^p}<\infty$ for some $p>2$.
\end{enumerate}
More precisely, (i) and (ii) are equivalent for fixed $\b>0$, (i) implies (iii)
for $2<p<\frac 52$ with $\frac 1p>\frac {2-\b}4$, and (iii) implies (i) with
$\b=1-\frac 2p$.

Moreover, (i) implies for some constant $C'$ on $D^*\times\S$ and
$B^*\times\S$ respectively
\begin{itemize}
\item[(iv)]
$\|F_\X(r,\p)\|_{L^\infty(\S)}\leq C' r^{\b-2} (\sin\p)^{-2}$ 
for all $r\leq r_0$, $\p\in (0,\pi)$.
\item[(iv')]
$\|F_\X(r,\p)\|_{L^\infty(\S)}\leq C' r^{\b-2}$ 
for all $r\leq r_0$, $\p\in [0,2\pi]$.
\end{itemize}
\end{lem}

\begin{rmk} \label{rmk decay}
If (\ref{bvp b}) or (\ref{bvp i}) in the above lemma are not satisfied, then still
$(ii)\Rightarrow (i)$, $(iii)\Rightarrow (i)$, and 
$(ii) \& (iv)\Rightarrow (iii)$ or $(ii) \& (iv')\Rightarrow (iii)$ respectively.
\end{rmk}

We will first show how this lemma and the subsequent theorem imply our main results, and then
give all proofs.
The following removal of singularities assumes a control of the curvature as given
by lemma~\ref{lem isoper} and \ref{lem decay} for finite energy solutions of (\ref{bvp b}) 
or (\ref{bvp i}).

\begin{thm} \hspace{1mm} \label{thm rem} \\
\vspace{-5mm} 
\begin{enumerate}
\item 
Let $\X\in\cA(B^*\times\S)$ satisfy (ii) and (iv') of lemma~\ref{lem decay}
with some constant $\b>0$.
Assume in addition that there exists a gauge in which $\X=A+R\dr$ with $\P\equiv 0$.
Let $2<p<\frac 52$ such that $\frac 1p>\frac {2-\b}4$.
Then there exists $u\in\cG^{2,p}_{\rm loc}(B^*\times\S)$ such that $u^*\X$ extends to
a connection $\tX\in\cA^{1,p}(B\times\S)$.

Moreover, if $\X$ is anti-self-dual, then $\tX$ will also be anti-self-dual.
\item 
Let $\X\in\cA(D^*\times\S)$ satisfy (ii) and (iv) of lemma~\ref{lem decay}
with some constant $\b>0$.
Let $2<p<\frac 52$ such that $\frac 1p>\frac {2-\b}4$.
Then there is \hbox{$u\in\cG^{2,p}_{\rm loc}(D^*\times\S)$} such that $u^*\X$ extends to
a connection $\tX\in\cA^{1,p}(D\times\S)$.

Moreover, if $\X$ satisfies (\ref{bvp b}), then $\tX$ will be a solution of (\ref{bvp}).
\end{enumerate}
\end{thm}

\noindent
{\bf Proof of theorem \ref{thm B} and remark \ref{rmk C}: } \\
Let $\X\in\cA(D^*\times\S)$ satisfy (\ref{bvp b}) and have finite energy 
$\E(r_0)<\infty$. Then lemma~\ref{lem isoper}~(ii) implies that $\E(r)\leq C r^{2\b}$
with $\b>0$, and hence we also have (ii) and (iv) as in lemma~\ref{lem decay}.
Now pick any $2<p<\frac 52$, and in case $0<\b<2$ choose it such that
$p<\frac 4{2-\b}$. Then theorem~\ref{thm rem}~(ii) provides a 
gauge transformation $u\in\cG^{2,p}_{\rm loc}(D^*\times\S)$ such that 
$u^*\X=\tX|_{D^*\times\S}$, where $\tX\in\cA^{1,p}(D\times\S)$ is a solution of (\ref{bvp}).
By the regularity \cite[Thm~A]{W elliptic} for solutions of (\ref{bvp}) we can multiply
$u$ by another gauge transformation in \hbox{$\cG^{2,p}(D\times\S)$}
(hence still $u\in\cG^{2,p}_{\rm loc}(D^*\times\S)$) such that $\tX\in\cA(D\times\S)$ is 
smooth.

Since on $D^*\times\S$ both $\X$ and $\tX$ are smooth and $u^*\X=\tX$
(i.e.\ $\rd u = u\tX - \X u $) we also know that $u\in\cG(D^*\times\S)$ is smooth. 

The proof of remark~\ref{rmk C} is exactly the same. Here lemma~\ref{lem isoper}~(i)
and theorem~\ref{thm rem}~(i) require the assumption that $\X\in\cA(B^*\times\S)$
is gauge equivalent to a connection with $\P\equiv 0$.
Moreover, this argument only uses the wellknown regularity theorem for anti-self-dual 
connections (see e.g.\ \cite[Thm~9.4]{W}).
\QED

Lemma~\ref{lem decay} will be a consequence of the following mean value inequalities.

\begin{lem} \label{lem mean curv}
There exist constants $C$ and $\ep>0$ such that the following holds.
Let $\X$ be a smooth connection on $D^*\times\S$ or $B^*\times\S$ 
that satisfies (\ref{bvp b}) or (\ref{bvp i}) respectively.
Suppose that $\E(2r)\leq\ep$ for some $0<r\leq \half r_0$, then 
\begin{itemize}
\item[(i)] On $D^*\times\S$ and $B^*\times\S$
$\displaystyle\qquad\quad
\sup_\p \|F_\X(r,\p)\|_{L^2(\S)}^2 \leq C r^{-2} \E(2r) .
$
\item[(ii)] On $D^*\times\S$ for all $\p\in (0,\pi)$
$\displaystyle\qquad\;
\|F_\X(r,\p)\|_{L^\infty(\S)}^2 \leq C (r\sin\p)^{-4} \E(2r) .
$
\item[(ii')]On $B^*\times\S$ for all $\p\in[0,2\pi]$
$\displaystyle\qquad
\|F_\X(r,\p)\|_{L^\infty(\S)}^2 \leq C r^{-4} \E(2r) .
$
\end{itemize}
\end{lem}
\Pr
We prove (i) in three steps and deduce (ii) and (ii') in the fourth.

\medskip

\noindent{\bf Step 1:} {\it We find constants $C$ and $\ep>0$ such that under the
above assumptions }
$$
\sup_\p \, r \, \|F_\X(r,\p)\|_{L^2(\S)} \leq C .
$$
Assume that for some fixed $\ep>0$ (that we shall fix later on) there is no such bound $C$.
Then we find a sequence of smooth connections $\X^\n$ on $D^*\times\S$ or 
$B^*\times\S$ satisfying (\ref{bvp b}) or (\ref{bvp i}) respectively,
and we find $\br^\n \to r^\infty\in[0,\half r_0]$ and $\bph^\n\to\p^\infty$
such that $\E^\n(2\br^\n)\leq\ep$ but 
$\br^\n \, \|F_{\X^\n}(\br^\n,\bph^\n)\|_{L^2(\S)} \to \infty $.
Here $\E^\n(\cdot)$ denotes the energy function (\ref{energy b}) or 
(\ref{energy i}) of $\X^\n$.
Given this we can choose $0<\bep^\n\leq\half \br^\n$ such that $\bep^\n\to 0$ but still
$ \bep^\n \, \|F_{\X^\n}(\br^\n,\bph^\n)\|_{L^2(\S)} \to \infty $.
The Hofer trick, lemma~\ref{lem HT} then yields $0<\ep^\n\leq\bep^\n$ 
(in particular $\ep^\n\to 0$) 
and $(r^\n,\p^\n)\to (r^\infty,\p^\infty)$ such that the following holds:
Firstly, with $R^\n := 2 \| F_{\X^\n}(r^\n,\p^\n) \|_{L^2(\S)} \to \infty $
we have
$$
\ep^\n \, R^\n \;\geq\; 2 \bep^\n \, \|F_{\X^\n}(\br^\n,\bph^\n)\|_{L^2(\S)} \;\to\;\infty .
$$
Secondly,
$$
\bigl\| F_{\X^\n}(r,\p) \bigr\|_{L^2(\S)} 
\;\leq\; 2 \, \bigl\| F_{\X^\n}(r^\n,\p^\n) \bigr\|_{L^2(\S)} \;=\; R^\n 
\qquad\forall (r,\p)\in B_{\ep^\n}(r^\n,\p^\n) .
$$
Here $B_{\ep^\n}(r^\n,\p^\n)$ denotes the Euclidean ball, where just the center 
$(r^\n,\p^\n)$ is given in polar coordinates. 
It is contained in $B_{2\br^\n}^*$ because $|r^\n-\br^\n|\leq\bep^\n\leq\half\br^\n$.
Moreover, in the boundary case it is understood to be intersected with $D$, 
so it is contained in $D_{2\br^\n}^*$.
Now proposition~\ref{prp est}~(ii) (with a fixed metric and any $\D>0$) provides a constant $C$
such that for all sufficiently large $\n\in\N$
$$
\bigl\| F_{A^\n}(r,\p) \bigr\|_{L^\infty(\S)} \;\leq\; C ( R^\n )^2 
\qquad\forall (r,\p)\in B_{\frac 12 \ep^\n}(r^\n,\p^\n).
$$
Putting this into the estimate of lemma~\ref{lem Laplacians} we obtain 
on $B_{\frac 12 \ep^\n}(r^\n,\p^\n)$ 
$$
\laplace  \bigl\| F_{\X^\n} \bigr\|_{L^2(\S)}^2
\;\leq\; C ( R^\n )^2 \bigl\| F_{\X^\n} \bigr\|_{L^2(\S)}^2 
$$
with another constant $C$, and in the boundary case moreover
$$
-\tfrac\pd{\pd t}\bigr|_{t=0} \bigl\| F_{\X^\n} \bigr\|_{L^2(\Si)}^2 
\;\leq\; C \bigl( \bigl\| F_{\X^\n} \bigr\|_{L^2(\Si)}^2 + \bigl\| F_{\X^\n} \bigr\|_{L^2(\Si)}^3 \bigr) .
$$
Now we fix $\ep=\frac 13\m C^{-2}$ with the $\m>0$ from proposition~\ref{prp mean bdy}. 
Then due to $\E^\n(2\br^\n)\leq\ep$ the mean value inequality applies to the functions 
$\| F_{\X^\n}\|_{L^2(\Si)}^2$ and yields with a new constant $C'$
$$
\bigl\| F_{\X^\n} (r^\n,\p^\n) \bigr\|_{L^2(\Si)}^2 
\leq C' \bigl( (R^\n)^2 + (\ep^\n)^{-2} \bigr) 
\int_{B_{\frac 12 \ep^\n}(r^\n,\p^\n)} \bigl\| F_{\X^\n} \bigr\|_{L^2(\Si)}^2 .
$$
If we moreover choose $\ep\leq\frac 1{2C'}$, then this implies
$2 (R^\n)^2 \leq (R^\n)^2 + (\ep^\n)^{-2} $
and thus $ (\ep^\n R^\n)^2 \leq 1 $ in contradiction to $\ep^\n R^\n \to\infty$.

\medskip

\noindent{\bf Step 2:} {\it We find constants $C$ and $\ep>0$ such that under the
above assumptions }
$$
\sup_\p \, r^2 \, \|F_A(r,\p)\|_{L^\infty(\S)} \leq C .
$$
Again arguing by contradiction we find a sequence of smooth connections $\X^\n$ on $D^*\times\S$ or 
$B^*\times\S$ satisfying (\ref{bvp b}) or (\ref{bvp i}) respectively, moreover 
$r^\n \to r^\infty\in[0,\half r_0]$ and $\p^\n\to\p^\infty$
such that $\E^\n(2 r^\n)\leq\ep$ but $ (r^\n)^2 \, \|F_{A^\n}(r^\n,\p^\n)\|_{L^\infty(\S)} \to \infty $.

Let $0<\ep^\n\leq \frac 12 r^\n$, then we know from step~1 that for some $\D>0$
$$
\|F_{\X^\n}(r,\p)\|_{L^2(\S)} \leq  2 (r^\n)^{-1} \D 
\qquad\forall (r,\p)\in B_{\ep^\n}(r^\n,\p^\n) .
$$
Now choose $R^\n\geq 2 (r^\n)^{-1} \D$ such that $R^\n\to\infty$, then the above
is true with $\ep^\n=\D (R^\n)^{-1}\leq\half r^\n$. 
Furthermore,  $\ep^\n\to 0$ and $\ep^\n R^\n = \D > 0$. 
So proposition~\ref{prp est}~(ii) asserts that for sufficiently large $\n\in\N$ 
and some constant $C$
$$
\|F_{A^\n}(r^\n,\p^\n)\|_{L^\infty(\S)} 
\;\leq\;  C (R^\n)^2 
\;=\; 4 C \D^2  (r^\n)^{-2} 
$$
in contradiction to $ (r^\n)^2 \, \|F_{A^\n}(r^\n,\p^\n)\|_{L^\infty(\S)} \to \infty $.

\pagebreak

\noindent{\bf Step 3:} {\it Proof of (i)}\\
Fix a connection $\X$ as assumed and consider a point $(r,\p)$ with $\E(2r)\leq\ep$.
Here we first choose $\ep>0$ as in step~2.
The $L^\infty$-bound from step~2 can be put into the estimate of 
lemma~\ref{lem Laplacians} to find another constant $C$ such that on 
$B_{\frac 12 r}(r,\p)$ 
$$
\laplace  \bigl\| F_\X \bigr\|_{L^2(\S)}^2
\;\leq\; C  r^{-2} \bigl\| F_\X \bigr\|_{L^2(\S)}^2 .
$$
In the boundary case this lemma also provides
$$
-\tfrac\pd{\pd t}\bigr|_{t=0} \bigl\| F_\X \bigr\|_{L^2(\Si)}^2 
\;\leq\; C \bigl( \bigl\| F_\X \bigr\|_{L^2(\Si)}^2 + \bigl\| F_\X \bigr\|_{L^2(\Si)}^3 \bigr) .
$$
Now we can choose a smaller $\ep>0$ such that $\ep\leq \frac 13 \m C^{-2}$ with the $\m>0$ from 
proposition~\ref{prp mean bdy}. Then we obtain the following mean value inequality for the function 
$\| F_\X \|_{L^2(\Si)}^2$ with another constant $C'$,
$$
\bigl\| F_\X (r,\p) \bigr\|_{L^2(\Si)}^2 
\;\leq\; C' r^{-2} \int_{B_{\frac 12 r}(r,\p)} \bigl\| F_\X \bigr\|_{L^2(\Si)}^2 
\;\leq\; 2 C' r^{-2} \E(2r) .
$$
%
%
%
%
%
%
%
%

\medskip

\noindent{\bf Step 4:} {\it Proof of (ii),(ii')}\\
It suffices to prove the estimates for $r\leq \br_0$ with some fixed $\br_0>0$, 
since then in case $\br_0 < r \leq \half r_0$ (and similarly in the boundary case)
$$
\|F_\X(r,\p)\|_{L^\infty(\S)}^2 \;\leq\; C (\br_0)^{-4} \E(2\br_0) 
\;\leq\; C \Bigl(\frac{r_0}{2\br_0}\Bigr)^4  r^{-4} \E(2r) .
$$
First, let $\br_0>0$ be the minimum of the injectivity radius on $\S$ for the metrics
$g_{s,t}$. Then we choose $\br_0>0$ even smaller
such that the pullback of all these metrics under normal coordinates on a ball
of radius $\br_0$ is $\cC^1$-close to the \hbox{Euclidean} metric on $\R^2$.
Thus we will be able to work with uniform constants $C$ and $\m>0$ in
proposition~\ref{prp mean int}.

In the interior case we consider a connection $\X$ as assumed
and any point $(r,\p,z)\in B_{\br_0}^*\times\S$.
The normal coordinates centered at this point give a coordinate chart
on $B_{\frac 12 r}(0)\subset\R^4$.
From lemma~\ref{lem Laplacians} we have a uniform constant $C$ such that
on $B_{\frac 12 r}(r,\p,z)\subset B^*\times\S$ 
$$
\laplace \bigl| F_\X \bigr|^2
\;\leq\; C \bigl| F_\X \bigr|^2 + 8 \bigl| F_\X \bigr|^3 .
$$
Now let $0<\ep\leq\frac \m {65}$, then proposition~\ref{prp mean int} applies to the
pullback of the function $| F_\X |^2$ on the coordinate chart $B_{\frac 12 r}(0)$
and asserts that
$$
\bigl| F_\X (r,\p,z) \bigr|^2 
\;\leq\; C \bigl( 1 + r^{-4} \bigr) \int_{B_{\frac 12 r}(r,\p,z)} \bigl| F_\X \bigr|^2 
\;\leq\; C r^{-4} \E(2r) .
$$
Here $C$ denotes any finite constant and we have used $1 \leq (r_0)^4 r^{-4}$.

In the boundary case on $D^*\times\S$ we use the same mean value inequality
on the ball $B_\r(r,\p,z)\subset D^*\times\S$ of radius $\r=\half r \sin\p$
for any $(r,\p,z)\in D_{\br_0}^*\times\S$ and $0<\p<\pi$. 
The normal coordinates centered at $(r,\p,z)$ give a coordinate chart
on the full ball $B_\r(0)\subset\R^4$.
With the same estimate on $\laplace | F_\X |^2$ and the same
$\ep>0$ as above we then apply proposition~\ref{prp mean int} to obtain
$$
\bigl| F_\X (r,\p,z) \bigr|^2 
\;\leq\; C \bigl( 1 + \r^{-4} \bigr) \int_{B_\r(r,\p,z)} \bigl| F_\X \bigr|^2 
\;\leq\; C (r\sin\p)^{-4} \E(2r) .
$$
Again, $C$ denotes any finite constant, and $1 \leq (r_0)^4 (r\sin\p)^{-4}$.
\QED

\medskip

\noindent
{\bf Proof of lemma \ref{lem decay} and remark \ref{rmk decay}: }
We will use $C$ and $C'$ to denote all finite constants. 
These might depend on the connection $\X$.

\medskip

\noindent \underline{$(i)\Rightarrow (ii)$ :}\; 
Since $\E(r_0)<\infty$ we must have $\E(r)\to 0$ as $r\to 0$.
So we find $\br>0$ such that for all $0<r\leq\br$ we obtain from lemma~\ref{lem mean curv}
$$
\sup_\p \|F_\X(r,\p)\|_{L^2(\S)} \;\leq\; r^{-1} \sqrt{C\,\E(2r)}
 \;\leq\; C'r^{\b-1} .
$$
For $\br<r\leq r_0$ we have with a constant $C'$ depending on $\br$ or $r_0$ 
$$
\sup_\p \|F_\X(r,\p)\|_{L^2(\S)} 
\;\leq\; \sup_{r\in[\br,r_0]}\sup_\p \, \|F_\X(r,\p)\|_{L^2(\S)} 
\;=\; C \;\leq\; C' r^{\b-1} .
$$

\noindent \underline{$(ii)\Rightarrow (i)$ :}\; Without using (\ref{bvp b}) or (\ref{bvp i}) 
we can simply calculate for all $\r\leq\half r_0$
$$
\E(\r) 
\;=\; \half \int_0^\r \int \|F_\X(r,\p)\|_{L^2(\S)}^2  \,r\,\dph\,\dr 
\;\leq\; \pi \int_0^\r C^2 r^{2\b-1}  \,\dr 
\;\leq\; C' \r^{2\b} .
$$
This already implies $\E(r_0)<\infty$. Then for $\half r_0 <\r \leq r_0$ we have 
$$
\E(\r) \;\leq\; \E(r_0) \,(\half r_0)^{-2\b}  \r^{2\b} 
\;=\; C  \r^{2\b} .
$$

\noindent \underline{$(i)\Rightarrow (iv), (iv')$ :}\; 
Since $\X$ is smooth away from $\{0\}\times\S$ it suffices to establish the estimates
for all $0<r\leq\br$. We pick $\br>0$ such that the assumptions of
lemma~\ref{lem mean curv} are satisfied, in particular $\E(2\br)\leq\ep$.
Then in the boundary case and the interior case respectively the lemma asserts
\begin{align*}
\|F_\X(r,\p)\|_{L^\infty(\S)} &\;\leq\; (r\sin\p)^{-2} \sqrt{C\,\E(2r)} \;\leq\;
C' r^{\b-2} (\sin\p)^{-2}
\qquad\,\forall \p\in (0,\pi) , \\
\|F_\X(r,\p)\|_{L^\infty(\S)} &\;\leq\, \qquad r^{-2} \sqrt{C\,\E(2r)} \quad \;\leq\; C' r^{\b-2} 
\qquad\qquad \qquad\forall \p\in[0,2\pi] .
\end{align*}

\noindent \underline{$(i)\Rightarrow (iii)$ :}\; This works the same for $D^*\times\S$ and 
$B^*\times\S$, so we only consider the first case.
(In the second case, the $\sin\p$-factor can be dropped.)
We already know that (i) implies (ii) and (iv).
Then just working with these two assumptions, we can interpolate for all $p>2$
\begin{align*}
\bigl\| F_\X \bigr\|_{L^p(D^*\times\S)}^p
&\;=\; \lim_{\d\to 0} \int_{(D_{r_0}\setminus D_\d)\times\S} \bigl| F_\X \bigr|^p   \\
&\;=\; \lim_{\d\to 0} \int_\d^{r_0} \int_0^\pi \bigl\|F_\X (r,\p)\bigr\|_{L^\infty(\S)}^{p-2} 
                                               \bigl\|F_\X (r,\p)\bigr\|_{L^2(\S)}^2   \,r\,\dph\,\dr \\
&\;\leq\; \lim_{\d\to 0} \; C \int_\d^{r_0} \int_0^\pi r^{(\b-2)(p-2) + 2(\b-1)} 
                                                   (\sin\p)^{-2(p-2)} \,r\,\dph\,\dr \\
&\;\leq\; \lim_{\d\to 0} \; 2C  \int_0^{\frac\pi 2} \bigl(\tfrac 2\pi \p\bigr)^{-2(p-2)} \,\dph
                                \int_\d^{r_0}r^{(\b-2)p + 3} \,\dr  .
\end{align*}
Here we use $\sin\p\geq\frac 2\pi \p$ for $\p\in[0,\frac\pi 2]$.
The $\p$-integral is finite for $p<\frac 52$, and the $r$-integral converges to a finite value
if $(\b-2)p > -4$. So if $\b\geq 2$, then we just need $2<p<\frac 52$, and if $\b<2$, then
we need in addition $p<\frac 4{2-\b}$.

\medskip

\noindent \underline{$(iii)\Rightarrow (i)$ :}\; This is the same calculation for both 
$D^*\times\S$ and $B^*\times\S$, and it works without the assumption 
(\ref{bvp b}) or (\ref{bvp i}).
 In the first case for all $r\leq r_0$,
\begin{align*}
\E(r) 
&\;=\; \lim_{\d\to 0} \half  \int_{(D_r\setminus D_\d)\times\S} |F_\X|^2   \\
&\;\leq\; \half {\rm Vol}(D_r\times\S)^{1 - \frac 2p} \;
      \lim_{\d\to 0} \Bigl(\int_{(D_r\setminus D_\d)\times\S} |F_\X|^p \Bigr)^{\frac 2p}    
\;\leq\; C' r^{2 ( 1 -\frac 2 p )} .
\end{align*}

\vspace{-5.5mm}

\hspace*{\fill}$\Box$
\medskip

\noindent
{\bf Proof of theorem \ref{thm rem}: }
We will give the full proof in the boundary case (ii) and point out where it differs 
(mostly simplifies) in the interior case (i).

Given a connection $\X\in\cA(D^*\times\S)$ as assumed we first put it into the special
gauge $\X=A+R\dr$ with $A:D^*\to\cA(\S)$ and $R:D^*\to\cC^\infty(\S,\cg)$ such that
$R|_{\p=\frac \pi 2} \equiv 0$ (and $\P\equiv 0$).
This is achieved by a gauge transformation $u\in\cG(D^*\times\S)$ that is determined 
as follows: For every $z\in\S$ first solve
$\pd_r u = -R u$ with initial value $u(r_0,\frac \pi 2,z)=\one$,
to determine $u(\cdot,\frac \pi 2,z)$, then for each $r\in(0,r_0]$ 
use this as initial value and solve $\pd_\p u = -\P u$ to obtain
$u(r,\p,z)$ for all $\p\in[0,\pi]$. That way the gauge is fixed up to
a gauge transformation on $\S$, i.e.\ independent of $(r,\p)\in D^*$.
(In case (i) this construction does in general not yield $u(r,0,z)= u(r,2\pi,z)$ and
hence define a gauge transformation on $B^*\times\S$. Thus the existence of this
gauge is an assumption in the theorem. Given this gauge, one then only needs to solve
$\pd_r u = -R u$ at $\p=\frac \pi 2$.)
In this gauge and splitting, the norm of the curvature is
$$
\bigl| F_\X \bigr|^2 = \bigl| F_A \bigr|^2 + \bigl| \pd_r A - \rd_A R \bigr|^2 
+ r^{-2}\bigl| \pd_\p R \bigr|^2 + r^{-2} \bigl| \pd_\p A \bigr|^2 .
$$
In particular, note that
$$
\bigl| \pd_\p\X \bigr|^2 = \bigl| \pd_\p R \bigr|^2 + \bigl| \pd_\p A \bigr|^2 
\leq r^2 \bigl| F_\X \bigr|^2 ,
\qquad
\bigl| \pd_r\X \bigr|_{\p=\frac \pi 2}^2
= \bigl| \pd_r A \bigr|_{\p=\frac \pi 2}^2 
\leq \bigl| F_\X \bigr|_{\p=\frac \pi 2}^2 .
$$
Next, we can combine the assumptions (ii) and (iv) as in lemma~\ref{lem decay} 
to obtain for any $q>2$ (in case (i) even without the $\sin\p$-term) 
\begin{align*} 
\bigl\| F_\X (r,\p) \bigr\|_{L^q(\S)}^q
\;\leq\; \bigl\| F_\X (r,\p) \bigr\|_{L^\infty(\S)}^{q-2} 
         \bigl\| F_\X (r,\p) \bigr\|_{L^2(\S)}^{2}        
\;\leq\; C \, r^{2-(2-\b)q} (\sin\p)^{4-2q} .  
\end{align*}
By integrating this over $D^*$ we recover (iii) of lemma~\ref{lem decay}:
If $2<p<\frac 52$ (in case (i) we only need $p>2$) and $\frac 1p > \frac {2-\b}4$
then
\begin{align} 
\bigl\| F_\X \bigr\|_{L^p(D_\r\times\S)}^p
\;\leq\; C \, \r^{4-(2-\b)p} \;\underset{\r\to 0}\longrightarrow\; 0 .  \label{Lp full} 
\end{align}
Moreover, for any $q>2$ we can read off for all $0<r\leq r_0$ and $0<\p<\pi$
\begin{align} 
\bigl\| \pd_\p\X (r,\p) \bigr\|_{L^q(\S)} &\leq C \, r^{\frac 2q +\b-1} (\sin\p)^{\frac 4q -2} , 
\nonumber\\
\bigl\| \pd_r\X (r,\tfrac \pi 2) \bigr\|_{L^q(\S)} &\leq C \, r^{\frac 2q +\b-2} . \label{Lp slice} 
\end{align}
Integrating the second estimate shows that for $\frac 1p > \frac {1-\b}2$ there exists a limit
$\X(r,\frac\pi 2)=A(r,\frac\pi 2)\to A_0 \in \cA^{0,p}(\S)$ as $r\to 0$. The first estimate then 
implies $\X(r,\cdot) \to A_0$ in $\cC^0([0,\pi],\cA^{0,p}(\S))$.
This motivates the following construction:

Fix a smooth cutoff function $h:[0,\infty)\to [0,1]$ 
with $h|_{[0,\ep]}\equiv 0$ and $h|_{[1-\ep,\infty)}\equiv 1$ for some $\ep>0$ 
and such that $|h'|\leq 2$.
Now for every $0<\r\leq \half r_0$ we set $A_\r:=\X(\r,\tfrac\pi 2)\in\cA(\S)$ 
and define $\X^\r\in\cA(D\times\S)$ by
$$
\X^\r (r,\p) \;:=\; A_\r + h(\tfrac r\r) \bigl( \X(r,\p) - A_\r \bigr) .
$$
Note that $\X^\r|_{D\setminus D_\r} = \X|_{D\setminus D_\r}$. 
We will find gauges for some sequence $\X^{\r_i}$, $\r_i\to 0$ such that these connections 
converge $W^{1,p}$-weakly. The limit will then be the extended connection 
$\tX\in\cA^{1,p}(D\times\S)$, and the gauge transformations will converge on $D^*\times\S$ to 
$u\in\cG_{\rm loc}^{2,p}$ such that $u^*\X = \tX|_{D^*\times\S}$.
This weak limit will be a consequence of Uhlenbeck's weak compactness theorem, 
so we have to control the curvatures
\begin{align*} 
F_{\X^\r} 
&=\; \rd A_\r \;+\; h(\tfrac r\r) \bigl( \rd \X - \rd A_\r \bigr)
   \;-\; \tfrac 1\r h'(\tfrac r\r) \bigl( \X - A_\r \bigr)\wedge\dr  \\
&\quad  +\; \half [A_\r\wedge A_\r] 
   \;+\; \half h(\tfrac r\r)^2 \bigl[(\X-A_\r)\wedge(\X-A_\r)\bigr]
   \;+\; h(\tfrac r\r) \bigl[A_\r \wedge (\X-A_\r)\bigr] \\
&=\; \bigl( 1 - h(\tfrac r\r) \bigr) F_{A_\r} 
   \;+\; h(\tfrac r\r) F_\X 
   \;-\; \tfrac 1\r h'(\tfrac r\r) \bigl( \X - A_\r \bigr)\wedge\dr \\
&\quad  \;+\; \half \bigl(h(\tfrac r\r)^2 - h(\tfrac r\r) \bigr) \bigl[(\X-A_\r)\wedge(\X-A_\r)\bigr] .
\end{align*}
From (\ref{Lp full}) we know that $F_\X\in L^p(D\times\S)$.
Now we shall see that $F_{\X^\r}\to F_\X $ in $L^p(D\times\S)$ as $\r\to 0$ :
\begin{align*} 
\bigl\| F_{\X^\r} - F_\X \bigr\|_{L^p(D\times\S)}
&\leq \bigl\| F_{A_\r} \bigr\|_{L^p(D_\r\times\S)}
   + \bigl\| F_\X \bigr\|_{L^p(D_\r\times\S)} \\
&\quad   + \tfrac 2\r \bigl\| \X - A_\r \bigr\|_{L^p(D_\r\times\S)}
   + \bigl\| \X-A_\r \bigr\|_{L^{2p}(D_\r\times\S)}^2 .
\end{align*}
The second term on the right hand side converges to zero by (\ref{Lp full}).
For the first term we use (\ref{Lp slice}) and recall that $p>2$ such that
$\frac 1p > \frac{2-\b}4$, so
\begin{align*} 
\bigl\| F_{A_\r} \bigr\|_{L^p(D_\r\times\S)}^p
\;=\; \int_{D_\r} \bigl\| F_\X (\r,\tfrac \pi 2) \bigr\|_{L^p(\S)}^p
\;\leq\; \half\pi C \, r^{4-(2-\b)p}
\;\underset{\r\to 0}\longrightarrow\; 0 . 
\end{align*}
To control the other two terms we first calculate for general $q>2$,
assuming $q\neq 4$, $\frac 2q + \b \neq 1$, and denoting all constants by $C$
\begin{align} 
&\bigl\| \X - A_\r \bigr\|_{L^q(D_\r\times\S)}^q  \nonumber\\
&= \int_{D_\r} \biggl\| \int_\r^r \pd_r\X(t,\tfrac\pi 2) \,\dt
                     + \int_{\frac\pi 2}^\p \pd_\p\X(r,\th) \,\dth \biggr\|_{L^q(\S)}^q 
\label{Lq full} \\
&\leq C \int_0^\r \int_0^{\frac\pi 2} \biggl( \int_r^\r \, t^{\frac 2q +\b-2} \,\dt
\,+ \int_\p^{\frac\pi 2} \, r^{\frac 2q+\b-1} (\sin\th)^{\frac 4q-2} \,\dth \biggr)^q\,r\,\dph\,\dr 
\nonumber\\
&\leq C \int_0^\r \biggl( r \r^{2-(1-\b)q} + r^{3-(1-\b)q} 
  + r^{3-(1-\b)q} \int_0^{\frac\pi 2}
    \bigl(1-\bigl(\tfrac 2\pi \p\bigr)^{\frac 4q-1}\bigr)^q \,\dph \biggr)\,\dr \nonumber\\
&\leq C \r^{4-(1-\b)q} . \nonumber
\end{align}
Here we have used the fact that $\sin\th\geq\frac 2\pi \th$ for $\th\in[0,\frac\pi 2]$.
The $\p$-integral then gives a finite value for $q<5$
and the $r$-integral converges for $\frac 1q > \frac{1-\b}4$. 
For $\frac 2q + \b = 1$ we have to deal differently with the $t$-integral in (\ref{Lq full}),
but still
\begin{align*} 
\int_0^\r \biggl( \int_r^\r \, t^{-1} \,\dt \biggr)^q\,r\,\dr 
\;=\; \int_0^\r r \, \ln\bigl(\tfrac\r r\bigr)^q \,\dr 
\;=\; \int_1^\infty \r^2 e^{-2y} y^q \,\dy 
\;=\; C \r^2 . 
\end{align*}
So (\ref{Lq full}) holds for $2<q<5$ if $q\neq 4$ and $\frac 1q > \frac{1-\b}4$. 
These conditions are all satisfied for $q=p$ since $\frac{2-\b}4 > \frac{1-\b}2$.
So (\ref{Lq full}) implies
$$
\tfrac 2\r \bigl\| \X - A_\r \bigr\|_{L^p(D_\r\times\S)}
\;\leq\; C \r^{\frac 4p +\b - 2} 
\;\underset{\r\to 0}\longrightarrow\; 0 . 
$$
Finally, we can choose $q=2p$ in (\ref{Lq full}) since then $4<q<5$ and
$\frac 1q > \frac{2-\b}8 > \frac{1-\b}4$. If we also note that
$\frac 2p > \frac{2-\b}2 > 1-\b$, then this gives
$$
\bigl\| \X-A_\r \bigr\|_{L^{2p}(D_\r\times\S)}
\;\leq\; C \r^{\frac 2p+\b-1} 
\;\underset{\r\to 0}\longrightarrow\; 0 . 
$$
Thus we have checked that $\| F_{\X^\r} - F_\X \|_{L^p(D\times\S)} \to 0$ as $\r\to 0$,
and hence $\| F_{\X^\r} \|_{L^p(D\times\S)}$ must be bounded for $\r\in (0,\half r_0]$.
In order to apply Uhlenbeck's weak compactness theorem (\cite[Thm~1.5]{U2} or \cite[Thm~A]{W}),
we choose a closed subset $D_{\frac 12 r_0}\subset U \subset {\rm int}(D)$ with smooth boundary, 
and we denote $U^*=U\setminus\{0\}$.
Then for some sequence $\r_i\to 0$ there exist gauge transformations \hbox{$u_i\in\cG^{2,p}(U\times\S)$}
such that the gauge transformed connections $u_i^*\X^{\r_i}$ converge $W^{1,p}$-weakly
to some $\tX \in\cA^{1,p}(U\times\S)$.
On every compact subset $K\subset U^*\times\S$ we have $\|\X^{\r_i} - \X\|_{W^{1,p}(K)}\to 0$.
In particular both $\|\X^{\r_i} \|_{W^{1,p}(K)}$ and
$\|u_i^*\X^{\r_i} \|_{W^{1,p}(K)}$ are bounded and thus 
$\|u_i^{-1}\rd u_i  \|_{W^{1,p}(K)}$ is bounded.
Hence for some further subsequence, $u_i|_{U^*\times\S}$ converges to some 
$u\in\cG^{2,p}_{\rm loc}(U^*\times\S)$ in the $\cC^0$-topology and in the
weak $W^{2,p}$-topology on every compact subset (see e.g.\ \cite[Lemma~A.8]{W}).
Furthermore, $u^*\X|_{U^*\times\S}=\tX|_{U^*\times\S}$ since on every compact subset both are the weak 
$W^{1,p}$-limit of $u_i^*\X^{\r_i}$.

On $(D\setminus U)\times\S$ we can now choose an extension of $u$ and define $\tX=u^*\X$
to obtain the claimed gauge transformation $u\in\cG^{2,p}_{\rm loc}(D^*\times\S)$
and extension $\tX\in\cA^{1,p}(D\times\S)$ with $u^*\X=\tX|_{D^*\times\S}$
The interior case (i) is proven exactly the same way. 
Just the estimates are simplified due to the absence of the $\sin\p$-term.

Furthermore, if $\X$ is anti-self-dual, then in both cases we also know that
$\tX$ is anti-self-dual since 
$\| F_\tX + * F_\tX \|_{L^p(D\times\S)}
=\|F_{u^*\X} + * F_{u^*\X}  \|_{L^p(D\times\S)}=0$.
Finally, suppose that $\X$ has Lagrangian boundary values $\X|_{(s,0)\times\S}\in\cL_Y$
for all $0<|s|\leq r_0$. Since $\cL_Y$ is gauge invariant 
and $\X^\r|_{\{r\geq\r\}}=\X|_{\{r\geq\r\}}$ we thus know for every $0<|s|\leq r_0$ that
$u_i^*\X^{\r_i}|_{(s,0)\times\S}\in\cL_Y$ for all sufficiently large $i\in\N$.
Moreover, $u_i^*\X^{\r_i}$ is bounded in $W^{1,p}(D\times\S)$, and the embedding
$W^{1,p}(D\times\S)\hookrightarrow \cC^0(D,L^p(\S))$ is compact (see \cite[Lemma~2.5]{W elliptic}).
So some subsequence of $u_i^*\X^{\r_i}|_{(s,0)\times\S}$ converges in $\cA^{0,p}(\S)$
for all $-r_0\leq s \leq r_0$.
Since $\cL_Y\subset\cA^{0,p}(\S)$ is closed this implies
$\tX|_{(s,0)\times\S}\in\cL_Y$ for all $0<|s|\leq r_0$.
This also holds at $s=0$ since $\tX|_{(s,0)\times\S}\in\cA^{0,p}(\S)$ is a continuous path for 
$s\in[-r_0,r_0]$ by the embedding $W^{1,p}(D\times\S)\hookrightarrow\cC^0(D,L^p(\S))$.
\QED

\pagebreak

 \bibliographystyle{alpha}

\end{document}